\documentclass[a4paper]{article}
\usepackage{amsmath}\usepackage{epsf,amsfonts,amsthm}
\usepackage{fontenc,indentfirst, delarray,amsfonts,amsmath,amssymb}

\newcommand{\be}{\begin{equation}}
\newcommand{\ee}{\end{equation}}
\newcommand{\bea}{\begin{eqnarray*}}
\newcommand{\eea}{\end{eqnarray*}}

\newcommand{\p}{\partial}
\newcommand{\la}{\langle}
\newcommand{\ra}{\rangle}
\newcommand{\raa}{\rightarrow}
\newcommand{\Ci}{C^{\infty}}
\newcommand{\E}{\ell}
\newcommand{\N}{\mathbb{N}}

\newcommand{\R}{\mathbb{R}}

\newcommand{\n}{\nabla}

\newcommand{\lp}{\left(}
\newcommand{\rp}{\right)}

\newcommand{\op}[1]{\!\!\mathop{\rm ~#1}\nolimits}

\mathchardef\za="710B  %\alpha
\mathchardef\zb="710C  %\beta
\mathchardef\zg="710D  %\gamma
\mathchardef\zd="710E  %\delta
\mathchardef\zve="710F %\epsilon
\mathchardef\zz="7110  %\zeta
\mathchardef\zh="7111  %\eta
\mathchardef\zy="7112 %\theta
\mathchardef\zi="7113  %\iota
\mathchardef\zk="7114  %\kappa
\mathchardef\zl="7115  %\lambda
\mathchardef\zm="7116  %\mu
\mathchardef\zn="7117  %\nu
\mathchardef\zx="7118  %\xi
\mathchardef\zp="7119  %\pi
\mathchardef\zr="711A  %\rho
\mathchardef\zs="711B  %\sigma
\mathchardef\zt="711C  %\tau
\mathchardef\zu="711D  %\upsilon
\mathchardef\zf="711E %\phi
\mathchardef\zq="711F  %\chi
\mathchardef\zc="7120  %\psi
\mathchardef\zw="7121  %\omega
\mathchardef\ze="7122  %\varepsilon
\mathchardef\zvy="7123  %\vartheta
\mathchardef\zvw="7124  %\varomega
\mathchardef\zvr="7125 %\varrho
\mathchardef\zvs="7126 %\varsigma
\mathchardef\zvf="7127  %\varphi
\mathchardef\zG="7000  %\Gamma
\mathchardef\zD="7001  %\Delta
\mathchardef\zY="7002  %\Theta
\mathchardef\zL="7003  %\Lambda
\mathchardef\zX="7004  %\Xi
\mathchardef\zP="7005  %\Pi
\mathchardef\zS="7006  %\Sigma
\mathchardef\zU="7007  %\Upsilon
\mathchardef\zF="7008  %\Phi
\mathchardef\zW="700A  %\Omega

\newtheorem{theo}{Theorem}
\newtheorem{prop}{Proposition}

\newtheorem{defi}{Definition}

\begin{document}

\title{A First Approximation for Quantization of Singular Spaces}
\author{Norbert Poncin\footnote{University of Luxembourg, Campus Limpertsberg, Institute
of Mathematics, 162A, avenue de la Fa\"iencerie, L-1511 Luxembourg
City, Grand-Duchy of Luxembourg, E-mail: norbert.poncin@uni.lu.
The research of N. Poncin was supported by grant R1F105L10. This
author also thanks the Erwin Schr\"odinger Institute in Vienna for
hospitality and support during his visits in 2006 and 2007.} ,
Fabian Radoux\footnote{University of Luxembourg, Campus
Limpertsberg, Institute of Mathematics, 162A, avenue de la
Fa\"iencerie, L-1511 Luxembourg City, Grand-Duchy of Luxembourg,
E-mail: fabian.radoux@uni.lu. F. Radoux is indebted to the
Luxembourg Ministry of Culture, Higher Education and Research, for
grant BFR 06/077}, Robert Wolak\footnote{Jagiellonian University,
ulica Reymonta 4 30-059 Krakow, Poland, E-mail:
Robert.Wolak@im.uj.edu.pl.}}\maketitle

\begin{abstract}
Many mathematical models of physical phenomena that have been
proposed in recent years require more general spaces than
manifolds. When taking into account the symmetry group of the
model, we get a reduced model on the (singular) orbit space of the
symmetry group action. We investigate quantization of singular
spaces obtained as leaf closure spaces of regular Riemannian
foliations on compact manifolds. These contain the orbit spaces of
compact group actions and orbifolds. Our method uses foliation
theory as a desingularization technique for such singular spaces.
A quantization procedure on the orbit space of the symmetry group
- that commutes with reduction - can be obtained from
constructions which combine different geometries associated with
foliations and new techniques originated in Equivariant
Quantization. The present paper contains the first of two steps
needed to achieve these just detailed goals.
\end{abstract}



\section{Introduction}

Quantization of singular spaces is an emerging issue that has been
addressed in an increasing number of recent works, see e.g.
\cite{BHPSingReduDefoQuan}, \cite{JHQuantReduc},
\cite{JHSingQuant}, \cite{JRSQuantStrat},
\cite{HLSingQuanCommRedu}, \cite{PflDefoQuanSympOrbi} ...\\

One of the reasons for this growing popularity originates from
current developments in Theoretical Physics related with reduction
of the number of degrees of freedom of a dynamical system with
symmetries. Explicitly, if a symmetry Lie group acts on the phase
space or the configuration space of a general mechanical system,
the quotient space is usually a singular space, an orbifold or a
stratified space ... The challenge consists in the quest for a
quantization procedure for these singular spaces that in addition commutes with reduction.\\

In this work, we investigate quantization of singular spaces
obtained as leaf closure spaces of regular Riemannian foliations
of compact manifolds. These contain the orbit spaces of compact
group actions (see \cite{RichTransGeom}). We build a quantization
that commutes by
construction with projection onto the quotient.\\

Our method uses the foliation as desingularization of the orbit
space $M/{\bar{\cal F}}$, where $\bar{\cal F}$ is the singular
Riemannian foliation made up by the closures of the leaves of the
regular Riemannian foliation ${\cal F}$ on manifold $M$. More
precisely, we combine Foliation Theory with recent techniques from
Natural and Equivariant Quantization. Close match can indeed be
expected, as both topics are tightly connected with natural
bundles and natural operators.\\

Equivariant quantization, in the sense of C. Duval, P. Lecomte,
and V. Ovsienko, developed as from 1996, see \cite{LMT},
\cite{LO}, \cite{DLO}, \cite{PL}, \cite{BM}, \cite{DO},
\cite{BHMP}, \cite{BM2}. This procedure requires equivariance of
the quantization map with respect to the action of a
finite-dimensional Lie subgroup of the symmetry group
$\op{Diff}(\R^n)$ of configuration space $\R^n$. Equivariant
quantization has first been studied in Euclidean space, mainly for
the projective and conformal subgroups, then extended in 2001 to
arbitrary manifolds, see \cite{PL2}. An equivariant, or better, a
natural quantization on a smooth manifold $M$ is a vector space
isomorphism $$Q[\n]:\op{Pol}(T^*M)\ni s\to Q[\n](s)\in{\cal
D}(M)$$ that verifies some normalization condition and maps, in
this paper, a smooth function $s\in \op{Pol}(T^*M)$ of ``phase
space'' $T^*M$, which is polynomial along the fibers, to a
differential operator $Q[\n](s)\in{\cal D}(M)$ that acts on
functions $f\in\Ci(M)$ of ``configuration space'' $M$. The
quantization map $Q[\n]$ depends on the projective class $[\n]$ of
an arbitrary torsionless covariant derivative $\n$ on $M$, and it
is natural with respect to all its arguments and for the action of
the group $\op{Diff}(M)$ of all local diffeomorphisms of $M$, i.e.
$$Q[\zf^*\n](\zf^*s)(\zf^*f)=\zf^*\lp Q[\n](s)(f)\rp,$$ $\forall
s\in\op{Pol}(T^*M),\forall f\in\Ci(M),\forall\zf\in\op{Diff}(M).$
Existence of such natural and projectively invariant quantizations
has been investigated in several works, see e.g. \cite{MB}, \cite{MR}, \cite{SH}.\\

In Foliation Theory, one distinguishes different geometries
associated with a foliated manifold $(M,{\cal F})$ (defined by a
H$\ae$fliger cocycle), namely adapted geometry, foliated geometry,
and transverse geometry. We denote in this introduction objects of
the adapted (resp. foliated, transverse) ``world'' by $O_3$ (resp.
$O_2$, $O_1$), whereas objects of leaf closure space $M/{\bar{\cal
F}}$ are denoted by $O_0$. Ideally, geometric structures of level
$i$ project onto geometric structures of level $i-1$, so that
$p(O_i)=O_{i-1}$, if we agree to denote temporarily any of these
projections by $p$. Let us also recall that, roughly, adapted
objects are objects on $M$ with some special properties, foliated
objects are locally constant along the leaves and live in the
normal bundle of the foliation, and that transverse objects are
objects on the transverse manifold $N$, which are ${\cal
H}$-invariant, where transverse manifold $N$ and the holonomy
pseudo-group ${\cal H}$ depend on the chosen defining cocycle of
foliation ${\cal F}$. In order to build a quantization $Q_0$ on
$M/{\bar{\cal F}}$, which commutes with the projection onto this
singular space, we construct adapted, foliated, and transverse
quantizations $Q_3$, $Q_2$, and $Q_1$, in such a way that \be
Q_{i-1}[p\,\n_{i}](p\,s_i)(p\,f_i)=p\lp
Q_i[\n_i](s_i)(f_i)\rp,\quad \forall
i\in\{1,2,3\}.\label{QuanCommProj_i}\ee Hence,
$$Q_0[\n_0](s_0)(f_0)=Q_0[p^3\n_3](p^3s_3)(p^3f_3)=p^3\lp
Q_{3}[\n_{3}](s_3)(f_3)\rp.$$ Observe that adapted quantization
$Q_3$ quantizes objects on $M$, whereas singular quantization
$Q_0$ only quantizes the objects of $M/\bar{\cal F}$. Eventually,
quantization actually commutes with projection onto the quotient.\\

The proofs of the three stages mentioned in Equation
(\ref{QuanCommProj_i}) are not equally hard. Since foliated
geometric objects on a foliated manifold $(M,{\cal F})$ are in
$1$-to-$1$ correspondence with ${\cal H}$-invariant geometric
objects on the transverse manifold $N$ associated with the chosen
cocycle, it is clear that stage $Q_2$ -- $Q_1$ is quite obvious.
The passages $Q_3$ -- $Q_2$ between the ``big'' adapted and
``small'' foliated quantizations, as well as transition $Q_1$ --
$Q_0$ from transverse quantization to singular quantization are
much more intricate.\\

The present paper should be accessible for readers who are not
necessarily experts in both fields, Natural Quantization and
Foliation Theory. In order to limit the length of the article, we
publish the stages $Q_3$ -- $Q_2$ and $Q_1$ -- $Q_0$ in two
different works. This publication deals with the {\sl first
approximation $Q_3$ -- $Q_2$ for quantization of singular spaces}.

\section{Natural and projectively invariant
quantization}

The constructions of $Q_3,$ $Q_2$, and $Q_1$ are nontrivial
extensions to the adapted, foliated, and transverse contexts, of
the proof of existence of natural and projectively invariant
quantization maps on an arbitrary smooth manifold, see \cite{MR}.
In the present section, we concisely describe the basic ideas of this technique.\\

In the theory of star-products, A. Lichnerowicz extensively used
the standard ordering prescription $Q_{\op{aff}}(\n)$ associated
with a covariant derivative $\n$. More precisely, consider the
space ${\cal D}^k(\zG(E),\zG(F))$ of $k$th order differential
operators between the spaces of sections $\zG(E)$ and $\zG(F)$ of
two vector bundles $E,F\raa M$ over a manifold $M$ as well as the
corresponding symbol space $\zG({\cal S}^kTM\otimes E^*\otimes
F)$. If $\n$ is a covariant derivative on $E$, denote by
$\n^k:\zG(E)\ni f\raa \n^kf\in\zG({\cal S}^{k}T^*M\otimes E)$ the
iterated symmetrized covariant derivative. Normal ordering map
$Q_{\op{aff}}(\n)$ then associates to any symbol $s\in\zG({\cal
S}^kTM\otimes E^*\otimes F)$ a differential operator
$Q_{\op{aff}}(\n)(s)\in{\cal D}^k(\zG(E),\zG(F))$ defined on any
section
$f\in\zG(E)$ by \be Q_{\op{aff}}(\n)(s)(f):=i_s(\n^kf)\in\zG(F).\label{AffQuant}\ee\\

The following example allows understanding the idea, due to M.
Bordemann, see \cite{MB}, underlying the construction of natural
and projectively invariant quantizations $Q$ on a manifold $M$,
see above. Set $M=S^n$, where $S^n$ is the $n$-dimensional sphere,
and $G=\op{GL}(n+1,\R)$. The elements $g\in G$ act on $\R^{n+1}$,
$g:\R^{n+1}\ni x\to gx\in\R^{n+1}$, and on $S^n$, $\zf_g: S^n\ni
x\raa gx/\vert\vert gx\vert\vert\in S^n$, where notations are
self-explaining. Observe that
$\tilde{M}:=\R^{n+1}\backslash\{0\}\raa S^n=M$ is a bundle with
typical fiber $\R^+_0$, and note that all $g$ preserve the
canonical connection of $\R^{n+1}$, but that the induced $\zf_g$
do usually not preserve the canonical Levi-Civita connection on
$S^n$. It seems therefore natural to lift the complex situation on
$M$ to the simpler situation on $\tilde{M}$. Thus, in order to
define $Q[\n](s)(f)$, where, see above, $\n$ denotes a torsionless
covariant derivative on $M$, $s$ a symbol in $\zG({\cal
S}\,TM)\simeq\op{Pol}(T^*M)$, and $f$ a function in $\Ci(M)$, one
constructs natural and projectively invariant lifts
\be\n\raa\tilde{\n},\;s\raa\tilde{s},\;f\raa\tilde{f},\label{LiftsAll}\ee
then sets \be\lp
Q[\n](s)(f)\rp^{\tilde{}}:=Q_{\op{aff}}(\tilde{\n})(\tilde{s})(\tilde{f}),\label{NatuInvaQuanSper}\ee
where $Q_{\op{aff}}$ is the standard ordering. The point is that
the normal ordering prescription, see its definition, is
natural---but of course not projectively invariant---and that we
require naturality and projective invariance for all the lifts. It
immediately follows that $Q$ inherits these properties (if the
projection onto the base behaves properly).\\

One of the proofs of existence of natural and projectively
invariant quantizations on an arbitrary smooth manifold $M$ is
based on the preceding example $M=S^n$ and consists of four
stages. In order to ensure readability of this paper, we recall
some concepts that are basic for further investigations and we
briefly depict the mentioned four stages.

\subsection{Basic concepts}

Let $M'$ and $M''$ be two smooth manifolds, let $m'$ be a point in
$M'$, and $U'$ a neighborhood of $m'$. Two smooth functions
$f:U'\to M''$ and $g:U'\to M''$ have at $m'$ a contact of order
$\ge r$, $r\in\N$, if and only if $f(m')=g(m')=:m''$ and, for any
chart of $M'$ around $m'$ and any chart of $M''$ around $m''$, the
components of the local forms $F$ of $f$ and $G$ of $g$ have the
same partial derivatives up to order $r$ at $m'$. It is well-known
that it suffices that this condition be satisfied for one pair of
charts. The classes of equivalence relation ``contact of order
$\ge r$ at $m'$'' are the $r$-jets at $m'$.

Clearly, if we denote the coordinates of $M'$ around $m'$ by $Z$,
the $r$-jet $j^r_{m'}(f)$ at $m'$ of a function $f$ is
characterized by the package $(\p_Z^{\za}F^{i})(Z(m'))$,
$\vert\za\vert\le r,$ $i\in\{1,\ldots,n''\}$, $n''=\op{dim}M''.$
Of course, a change of coordinates entails a change of the
characterizing package of derivatives. If, for instance, if we
exchange current coordinates $X$ in target manifold $M''$ for new
coordinates $Y$, the current and new local forms $F(Z)=X(f(Z))$
and $F'(Z)=Y(f(Z))$ are related by $F'(Z)=Y(X(f(Z)))$, where, in
order to simply, we used notations from Physics. It follows that
\be
\p_{Z^{j}}F'^{i}=\p_{X^a}Y^{i}\,\p_{Z^j}F^a\label{JetCharChangeCoord1}\ee
and that
\be\p_{Z^jZ^k}F'^{i}=\p_{X^{a}X^b}Y^{i}\,\p_{Z^k}F^b\,\p_{Z^j}F^{a}+\p_{X^a}Y^{i}\,\p_{Z^jZ^k}F^a.\label{JetCharChangeCoord2}\ee
These formul\ae$\mbox{}$ will be needed below. Let us also recall
that, for fixed charts, the characterizing package of derivatives
of the jet $j^r_{m'}(h\circ f)$ of a compound map, is obtained,
roughly spoken, by composition of the limited Taylor expansions of
the local forms of $h$ and $f$, if one agrees to suppress the
terms that have order $>r$.\\

We denote by $P^r$, $r\in\N$, the natural functor of order $r$---
between the category of $n$-dimensional smooth manifolds $M$ and
immersions $\zf:M\to M'$ (or, equivalently, globally defined local
diffeomorphisms) and the category of fiber bundles and bundle
maps---the objects of which are the $r$th order frame bundles
$P^rM=\{j^r_0(f) |\, f:0\in U\subset \R^n\to M,
T_0f\in\op{Isom}(\R^n,T_{f(0)}M)\},$ and the morphisms of which
are the principal bundle morphisms $P^r\phi:P^rM\to P^rM'$ defined
by $(P^r\phi)(j^r_0(f))=j^r_0(\phi\circ f).$ The structure group
of principal bundle $P^rM$ is $G^r_n=\{j^r_0(\zvf)|\, \zvf:0\in
U\subset \R^n\to\R^n, \zvf(0)=0, T_0\zvf\in\op{GL}(n,\R)\}$ and
its action on $P^rM$, $j^r_0(f).j^r_0(\zvf):=j^r_0(f\circ\zvf)$,
is well-defined in view of the above remark on jets of compound
maps. Note that structure group $G^1_n$ of the principal bundle of
linear frames $P^1M=:LM$ is $G^1_n\simeq\op{GL}(n,\R)$. Remark
further that if $j^2_0(\zvf)\in G^2_n$ is characterized by
$(0^{i},A^{i}_k,S^{i}_{kl})$ and $j^2_0(f)\in P^2M$ is
characterized in coordinates $X$ around $m:=f(0)$ by
$(X^{i}(m),B^{i}_k,T^{i}_{kl})_X$, then
$j^2_0(f).j^2_0(\zvf)=j^2_0(f\circ\zvf)$ is characterized by \be
(X^{i}(m),B^{i}_k,T^{i}_{kl})_X\cdot(0^{i},A^{i}_k,S^{i}_{kl})=(X^{i}(m),B^{i}_aA^{a}_k,B^{i}_aS^{a}_{kl}+T^{i}_{ab}A^a_kA^b_l)_X\label{ActionChar}\ee

It is easily verified that the isotropy subgroup of
$[e_{n+1}]:=[(0,\ldots,0,1)^{\tilde{}}\,]\in\R P^n$ for the
canonical action of the projective group
$$\op{PGL}(n+1,\R)=\left\{\left(\begin{array}{cc}{\cal A} & h \\
\za & a\end{array}\right):{\cal A}\in\op{GL}(n,\R),\za\in\R^{n*},
h\in\R^n, a\in\R_0\right\}\slash\R_0\op{id}$$ on the
$n$-dimensional real projective space $\R P^n$, is $$H(n+1,\R)=\left\{\left(\begin{array}{cc}{\cal A} & 0 \\
\za & a\end{array}\right):{\cal A}\in\op{GL}(n,\R),\za\in\R^{n*},
a\in\R_0\right\}\slash\R_0\op{id},$$ and that $H(n+1,\R)$ acts
locally on $\R^n$ by affine fractional transformations that
preserve the origin. Hence, $H(n+1,\R)$ can be viewed as Lie
subgroup of structure group $G^r_n$.

\begin{prop}\label{InclusionUsualContext} The natural inclusion $I:H(n+1,\R)\to G^2_n$ reads
$I:\left[\left(\begin{array}{cc}{\cal
A}&0\\\za&1\end{array}\right)\right]\mapsto (0,{\cal
A}^{i}_{j},-{\cal A}^{i}_j\za_k-{\cal A}^{i}_{k}\za_j)$.
\end{prop}

\begin{proof} The natural action of an element
$\left[\left(\begin{array}{cc}{\cal
A}&0\\\za&1\end{array}\right)\right]\in
H(n+1,\R)\subset\op{PGL}(n+1,\R)$ on $Z\in U\subset\R^n$, where
$U$ is a sufficiently small neighborhood of $0$, is $\frac{{\cal
A}Z}{\za Z+1}\in\R^n$. A short and easy computation then shows
that the second jet at $0$ of map $\zvf:Z\mapsto\frac{{\cal
A}Z}{\za Z+1}$ is characterized in canonical coordinates by $(0,
{\cal A}^{i}_j,-{\cal A}^{i}_j\za_k-{\cal A}^{i}_k\za_j)$.
\end{proof}

\subsection{Stage 1: Cartan bundle}

A {\sl projective structure} on a smooth manifold $M$ is a class
$[\n]$ of all torsion-free linear connections $\n'$ on $M$ that
are projectively equivalent to $\n$, i.e. that have the same
geometric geodesics as $\n$, or better still, that verify
\be\n'_XY-\n_XY=\za(X)Y+\za(Y)X,\label{ProjEquiv}\ee for all
$X,Y\in\op{Vect}(M)$ and some fixed $\za\in\zW^1(M)$ (H. Weyl).\\

The next theorem contains the first of two essential observations,
see \cite{MR}, that allow solving the problem of the
aforementioned projectively invariant lift $\n\to\tilde{\n}$, i.e.
that allow associating a unique connection to each projective
structure.

\begin{theo} Let $M$ be a smooth manifold. There is a canonical
$1$-to-$1$ correspondence between projective structures $[\n]$ on
$M$ and reductions $P=P(M,H(n+1,\R))$ to structure group
$H(n+1,\R)$ of the principal bundle $P^2M=P^2M(M,G^2_n)$ of second
order frames on $M$.\label{ProjStruReductions}\end{theo}

In the sequel, we refer to the bundles $P=P(M,H(n+1,\R))$ as {\sl
Cartan bundles}.

\subsection{Stage 2: Cartan connection}

The second observation then settles the question of connection
lift $[\n]\to\tilde{\n}$:

\begin{theo} A unique normal Cartan connection is associated with every
Cartan bundle $P(M,H(n+1,\R))$ of $M$.
\end{theo}

Let $G$ be a Lie group, $H$ a closed subgroup, $\frak{g}$ and
$\frak{h}$ the corresponding Lie algebras, and let $P=P(M,H)$
denote a principal $H$-bundle over a manifold $M$, such that
$\op{dim}M=\op{dim}G\slash H$. In this setting, a {\sl Cartan
connection} on $P(M,H)$ is a differential $1$-form
$\zw\in\zW^1(P)\otimes\frak{g}$ valued (not in Lie algebra ${\frak
h}$, but) in Lie algebra ${\frak g}$, which verifies the usual
requirements for connection $1$-forms, i.e. $${\frak
r}_s^*\zw=\op{Ad}(s^{-1})\zw\quad\mbox{and}\quad\zw(X^h)=h,$$
where ${\frak r_s}$ denotes the right action by $s\in H$ and where
$X^h$ is the fundamental vector field associated with $h\in
\frak{h}$. However, a third condition asks that
$\zw_u:T_uP\raa\frak{g}$ be a vector space isomorphism for any
$u\in P$. Hence, we have $\op{ker}\zw_u=0$, so that the basic
difference with Ehresmann connections is the absence of a
horizontal subbundle. For instance, if $H$ is a closed subgroup of
a Lie group $G$, the canonical Maurer-Cartan form is a Cartan
connection on the principal bundle $G(G\slash H,H)$.

\subsection{Stage 3: Lifts of symbols and functions}

In view of the preceding remarks, the role of connection lift
$\tilde{\n}$ to bundle $\tilde{M}$, see Equation (\ref{LiftsAll}),
is played by the unique Cartan connection $\zw$ associated with
the unique Cartan bundle $P=P(M,H)$, $H=H(n+1,\R)$, defined by the
considered projective structure $[\n]$ on $M$. Lifting symbols
$s\in\zG({\cal S}^kTM)$ and in particular functions $f\in\Ci(M)$
to objects $\tilde{s}$ and $\tilde{f}$ of $\tilde{M}\simeq P$, is
then quite obvious. Indeed, we have $\zG({\cal
S}^kTM)=\Ci(P^1M,{\cal S}^k\R^n)_{\op{GL}(n,\R)}$, where the
{\small RHS} denotes the space of $\op{GL}(n,\R)$-invariant ${\cal
S}^k\R^n$-valued functions of the linear frame bundle $P^1M$.
Since, there are canonical projections $P\subset P^2M\to P^1M$ and
$H\subset G^2_n\to G^1_n=\op{GL}(n,\R),$ it is easily seen that
\be\zG({\cal S}^kTM)=\Ci(P^1M,{\cal
S}^k\R^n)_{\op{GL}(n,\R)}\simeq \Ci(P,{\cal
S}^k\R^n)_H,\label{LiftArgu}\ee where the $H$-action on ${\cal
S}^k\R^n$ is induced by the corresponding $\op{GL}(n,\R)$-action.

\subsection{Stage 4: Construction of a natural and invariant quantization}

Equations (\ref{AffQuant}) and (\ref{NatuInvaQuanSper}) suggest
defining a natural and projectively invariant quantization on a
smooth manifold $M$, endowed with a projective structure
$[\n]\simeq P$, by \be\lp Q[\n](s)(f)\rp
^{\tilde{}}=Q_{\op{aff}}(\zw)(\tilde{s})(\tilde{f})=
i_{\tilde{s}}(\n^{\zw})^
k\tilde{f},\label{NatuProjInvaQuanFirsAtte}\ee where $\n^{\zw}$
denotes a covariant derivative associated with connection $1$-form
$\zw$. Whereas $\n^{\zw}$ can easily be defined, it turns out that
the {\small RHS} of Equation (\ref{NatuProjInvaQuanFirsAtte}) is a
function of $P$ that is not $H$-invariant, so that it does not
project onto a function $Q[\n](s)(f)$ of $M$, see Equation
(\ref{LiftArgu}), set $k=0$, and note that $\sim$ is just the
isomorphism $\simeq$. The solution consists in the substitution to
$\tilde{s}\in \Ci(P,{\cal S}^k\R^n)_H$ of a linear combination of
lower degree terms. These are obtained from tensor field
$\tilde{s}$ by means of a degree-lowering divergence operator
$\op{Div}^{\zw}=\sum_ji_{\ze^j}\n^{\zw}_{e_j}$, where $(e_j)_j$
and $(\ze^{j})_j$ are the canonical bases of $\R^n$ and $\R^{n *}$
respectively. Eventually, it can be proven, see \cite{MR}, that
$$\lp
Q[\n](s)(f)\rp^{\tilde{}}=\sum_{\E=0}^kc_{k\E}\;i_{\lp\op{Div}^{\zw}\rp^{\E}
\tilde{s}}\lp\n^{\zw}\rp^{k-\E} \tilde{f}$$ defines a natural and
projectively invariant quantization on $M$, if the coefficients
$c_{k\E}\in\R$ have some precise values.\\

In the following, we study extensions of the just detailed modus
operandi to the adapted and foliated geometries associated with
foliated manifolds.

\section{Adapted and foliated projective structures}

In this section, we investigate the link between adapted (resp.
foliated) projective structures and reductions of the principal
bundle of adapted (resp. foliated) second order frames.

\subsection{Adapted and foliated connections}

Let $(M,{\cal F})$ be a foliated manifold, more precisely, let $M$
be an $n$-dimensional smooth manifold endowed with a {\sl regular
foliation} ${\cal F}$ of dimension $p$ (and codimension $q=n-p$).
It is well-known that such a foliation can be defined as an
involutive subbundle $T{\cal F}\subset TM$ of constant rank $p$.

Foliation ${\cal F}$ can also be viewed as a partition into
(maximal integral) $p$-dimensional smooth submanifolds or {\sl
leaves}, such that in appropriate or {\sl adapted charts}
$(U_i,\zf_i)$ the connected components of the traces on $U_i$ of
these leaves lie in $M$ as $\R^p$ in $\R^n$ [pages of a book],
with transition diffeomorphisms of type $\psi_{ji}=\zf_{j}\circ
\zf_{i}^{-1}:\zf_i(U_{ij})\ni(x,y)\raa
(\psi_{ji,1}(x,y),\psi_{ji,2}(y))\in\zf_j(U_{ji})$,
$U_{ij}=U_i\cap U_j$ [the $\psi_{ji}$ map a page onto a page]. The
pages provide by transport to manifold $M$ the so-called {\sl
plaques} or {\sl slices} and these glue together from chart to
chart---in the way specified by the transition
diffeomorphisms---to give maximal connected injectively immersed
submanifolds, precisely the leaves of the foliation.

Eventually, foliation ${\cal F}$ can be described by means of a
H$\ae$fliger {\sl cocycle} ${\cal U}=(U_i,f_i,g_{ij})$ modelled on
a $q$-dimensional smooth manifold $N_0$. The $U_i$ form an open
cover of $M$ and the $f_i:U_i\raa f_i(U_i)=:N_i\subset N_0$ are
submersions that have connected fibers [the connected components
of the traces on the $U_i$ of the leaves of ${\cal F}$] and are
subject to the transition conditions $g_{ji}f_i=f_j$, where the
$g_{ji}:f_i(U_{ij})=:N_{ij}\raa N_{ji}:=f_j(U_{ji})$ are
diffeomorphisms that verify the usual cocycle condition
$g_{ij}g_{jk}=g_{ik}$. We refer to the disjoint union $N=\amalg_i
N_i$ as the (smooth, $q$-dimensional) {\sl transverse manifold}
and to ${\cal H}:=\ra g_{ij}\la$ as the pseudogroup of (locally
defined) diffeomorphisms or {\sl holonomy pseudogroup} associated
with the chosen cocycle ${\cal U}$.\\

A vector field $X\in\op{Vect}(M)$, such that $[X,Y]\in\zG(T{\cal
F})$, for all $Y\in\zG(T{\cal F})$, is said to be {\sl adapted}
(to the foliation). The space $\op{Vect}_{\cal F}(M)$ of adapted
vector fields is obviously a Lie subalgebra of the Lie algebra
$\op{Vect}(M)$, and the space $\zG(T{\cal F})$ of {\sl tangent}
(to the foliation) vector fields is an ideal of $\op{Vect}_{\cal
F}(M)$. The quotient algebra $\op{Vect}(M,{\cal
F})=\op{Vect}_{\cal F}(M)/\zG(T{\cal F})$ is the algebra of {\sl
foliated} vector fields.

Let $(x,y)$ be local coordinates of $M$ that are adapted to ${\cal
F}$, i.e. $x=(x^1,\ldots,x^p)$ are leaf coordinates and
$y=(y^1,\ldots,y^q)$ are transverse coordinates. The local form of
an arbitrary (resp. tangent, adapted, foliated) vector field is
then $X=\sum_{\zi=1}^p X^{\zi}(x,y)\p_{\zi}+\sum_{{\frak i}=1}^q
X^{{\frak i}}(x,y)\p_{\frak i}$, $\p_{\zi}=\p_{x^{\zi}}$,
$\p_{\frak i}=\p_{y^{\frak i}}$ (resp. $X=\sum_{\zi=1}^p
X^{\zi}(x,y)\p_{\zi}$, \be X=\sum_{\zi=1}^p
X^{\zi}(x,y)\p_{\zi}+\sum_{{\frak i}=1}^q X^{\frak i}(y)\p_{\frak
i},\label{AdapVectFiel}\ee \be[X]=[\sum_{{\frak i}=1}^q X^{\frak
i}(y)\p_{\frak i}],\label{FoliVectFiel}\ee where $[.]$ denotes the
classes in the aforementioned quotient algebra).\\

In Foliation Theory, vocabulary is by no means uniform. Let us
stress that adapted and foliated vector fields, see Equations
(\ref{AdapVectFiel}) and (\ref{FoliVectFiel}), may be viewed as
prototypes of all adapted and foliated structures used in this
paper.

For instance, a smooth function $f\in\Ci(M)$ is foliated (or
basic) if and only if $L_Yf=0,\forall Y\in\zG(T{\cal F}).$ We
denote by $\Ci(M,{\cal F})$ the space of all foliated functions of
$(M,{\cal F})$. A differential $k$-form $\zw\in\zW^k(M)$ is
foliated (or basic) if and only if
$\op{i}_Y\zw=\op{i}_Y\op{d}\zw=0,\forall Y\in\zG(T{\cal F}),$
where notations are self-explaining. Again, we denote by
$\zW^k(M,{\cal F})$ the space of all foliated differential
$k$-forms of $(M,{\cal F})$.

It is easily checked that $\Ci(M,{\cal F})\times\op{Vect}(M,{\cal
F})\ni (f,[X])\to f[X]:=[fX]\in\op{Vect}(M,{\cal F})$ defines a
$\Ci(M,{\cal F})$-module structure on $\op{Vect}(M,{\cal F})$.
Furthermore, $\op{Vect}(M,{\cal F})\times \Ci(M,{\cal F})\ni
([X],f)\to L_{[X]}f:=L_Xf\in\Ci(M,{\cal F})$ is the natural action
of foliated vector fields on foliated functions. Eventually, the
contraction of a foliated $1$-form $\za\in\zW^1(M,{\cal F})$ and a
foliated vector field $[X]\in \op{Vect}(M,{\cal F})$ is a foliated
function $\za([X]):=\za(X)\in\Ci(M,{\cal F})$.

\begin{defi} Let $(M,{\cal F})$ be a foliated manifold. An {\sl adapted connection}
$\n_{\cal F}$ is a linear torsion-free connection on $M$, such
that $\n_{\cal F}:\op{Vect}_{\cal F}(M)\times \zG(T{\cal F})\to
\zG(T{\cal F})$ and $\n_{\cal F}:\op{Vect}_{\cal
F}(M)\times\op{Vect}_{\cal F}(M)\to\op{Vect}_{\cal F}(M)$.
\end{defi}

\noindent {\bf Remark} In the following, we use the Einstein
summation convention, and, as already adumbrated above, Latin
indices $i,k,l\ldots$ (resp. Greek indices $\zi,\zk,\zl\ldots$,
German indices ${\frak i},{\frak k},{\frak l}\ldots$) are
systematically and implicitly assumed to vary in $\{1,\ldots,n\}$
(resp. $\{1,\ldots,p\}$, $\{1,\ldots,q\}$).\\

As torsionlessness means that $\n_{{\cal F},Y}X=\n_{{\cal
F},X}Y+[Y,X]$, it follows that $\n_{\cal F}:\zG(T{\cal
F})\times\op{Vect}_{\cal F}(M)\to\zG(T{\cal F})$.

Further, locally, in adapted coordinates, we have $\n_{{\cal
F},X}Y=\lp X^{i}\p_iY^k+\zG^k_{il}X^{i}Y^l\rp\p_k$, so that
condition $\n_{\cal F}:\op{Vect}_{\cal F}(M)\times \zG(T{\cal
F})\to \zG(T{\cal F})$ means that \be \zG^{\frak
k}_{i\zl}=\zG^{\frak k}_{\zl i}=0,\label{AdapConnLocal1}\ee
whereas condition $\n_{\cal F}:\op{Vect}_{\cal
F}(M)\times\op{Vect}_{\cal F}(M)\to\op{Vect}_{\cal F}(M)$ is then
automatically verified provided that Christoffel's symbols
$\zG^{\frak k}_{{\frak i}{\frak l}}$ are independent of $x$,
$\zG^{\frak k}_{{\frak i}{\frak l}}=\zG^{\frak k}_{{\frak i}{\frak
l}}(y)$.

\begin{defi} Consider a foliated manifold $(M,{\cal F})$. A {\sl
foliated} torsion-free connection $\n({\cal F})$ on $(M,{\cal F})$
is a bilinear map $\n({\cal F}):\op{Vect}(M,{\cal
F})\times\op{Vect}(M,{\cal F})\to\op{Vect}(M,{\cal F})$, such
that, for all $f\in\Ci(M,{\cal F})$ and all
$[X],[Y]\in\op{Vect}(M,{\cal F})$, the following conditions hold
true:
\begin{itemize}\item $\n({\cal F})_{f[X]}[Y]=f\n({\cal
F})_{[X]}[Y]$,\item $\n({\cal F})_{[X]}(f[Y])=\lp L_{[X]}f\rp
[Y]+f\n({\cal F})_{[X]}[Y],$\item $\n({\cal F})_{[X]}[Y]=\n({\cal
F})_{[Y]}[X]+[[X],[Y]].$
\end{itemize}
\end{defi}

In view of the above definitions, the local form (in adapted
coordinates (x,y)) of a foliated vector field is $[X]=X^{\frak
i}[\p_{\frak i}]$, $X^{\frak i}=X^{\frak i}(y)$, and a foliated
connection reads $$\n({\cal F})_{[X]}[Y]=X^{\frak i}\lp
L_{[\p_{\frak i}]}Y^{\frak l}\rp [\p_{\frak l}]+X^{\frak
i}Y^{\frak l}\;\zG({\cal F})^{\frak k}_{{\frak i}{\frak
l}}\;[\p_{\frak k}],\quad \zG({\cal F})^{\frak k}_{{\frak i}{\frak
l}}=\zG({\cal F})^{\frak k}_{{\frak i}{\frak l}}(y).$$

\begin{prop}\label{AdapProjEquiv} If two adapted connections $\n_{\cal F}$ and $\n'_{\cal F}$ of
a foliated manifold $(M,{\cal F})$ are projectively equivalent,
the corresponding differential $1$-form $\za\in\zW^1(M)$ is
foliated, i.e. $\za\in\zW^1(M,{\cal F})$.\end{prop}

\begin{proof} In adapted local coordinates $(x,y)$, projective equivalence of
$\n_{\cal F}$ and $\n'_{\cal F}$ reads
$(\zG'^k_{il}-\zG^k_{il})X^{i}Y^l=\za_iX^{i}Y^k+\za_iY^{i}X^k,\forall
k$. When writing this equation for $X^{i}=\zd^{i}_{\zi}$,
$Y^{l}=\zd^{l}_{\frak l}$, and $k={\frak l}$, we get, in view of
Equation (\ref{AdapConnLocal1}), $\za_{\zi}=0$. If we now choose
$X^{i}=\zd^{i}_{\frak i}$, $Y^{l}=\zd^{l}_{\frak l}$, and
$k={\frak i}\neq {\frak l}$, we finally see that $\za_{\frak l}$
is independent of $x$. \end{proof}

The following proposition is well-known:

\begin{defi} Two foliated connections $\n({\cal F})$ and $\n'({\cal F})$
of a foliated manifold $(M,{\cal F})$ are projectively equivalent,
if and only if there is a foliated $1$-form $\za\in\zW^1(M,{\cal
F})$, such that, for all $[X],[Y]\in\op{Vect}(M,{\cal F})$, one
has $\n'({\cal F})_{[X]}[Y]_-\n({\cal
F})_{[X]}[Y]=\za([X])[Y]+\za([Y])[X]$.\end{defi}

Eventually, adapted connections induce foliated connections.

\begin{prop}\label{AdapConnIndFolConn} Let $(M,{\cal F})$ be a foliated manifold of
codimension $q$. Any adapted connection $\n_{\cal F}$ of $M$
induces a foliated connection $\n({\cal F})$, defined by $\n({\cal
F})_{[X]}[Y]:=[\n_{{\cal F},X}Y]$. In adapted coordinates,
Christoffel's symbols $\zG({\cal F})^{\frak i}_{{\frak k}{\frak
l}}$ of $\n(\cal F)$ coincide with the corresponding Christoffel
symbols $\zG_{{\cal F},{\frak k}{\frak l}}^{\frak i}$ of $\n_{\cal
F}$. Eventually, projective classes of adapted connections induce
projective classes of foliated connections.\end{prop}

\begin{proof} It immediately follows from the definition of adapted
connections that for any $[X],[Y]\in\op{Vect}(M,{\cal F})$, the
class $\n({\cal F})_{[X]}[Y]:=[\n_{{\cal
F},X}Y]\in\op{Vect}(M,{\cal F})$ is well-defined. All properties
of foliated connections are obviously satisfied. If $(x,y)$ are
adapted coordinates, we have $\zG({\cal F})^{\frak k}_{{\frak
i}{\frak l}}[\p_{\frak k}]=\n({\cal F})_{[\p_{\frak i}]}[\p_{\frak
l}]=[\n_{\cal F,\p_{\frak i}}\p_{\frak l}]=[\zG_{\cal F,{\frak
i}{\frak l}}^{\frak k}\p_{\frak k}]=\zG_{\cal F,{\frak i}{\frak
l}}^{\frak k}[\p_{\frak k}]$, since $\zG_{{\cal F},{\frak i}{\frak
l}}^{\zk}=0$ and $\zG_{{\cal F},{\frak i}{\frak l}}^{\frak
k}=\zG_{{\cal F},{\frak i}{\frak l}}^{\frak k}(y)$. The remark on
projective structures follows immediately from preceding
observations.
\end{proof}

\subsection{Adapted and foliated frame bundles}

\subsubsection{Adapted frame bundles}

Since an adapted linear frame is a frame $(v_1,\ldots,v_{p+q})$ of
a fiber $T_mM$, $m\in M$, the first vectors $(v_1,\ldots,v_p)$ of
which form a frame of $T_m{\cal F}$, we denote by $P^r_{\cal F}M$
the principal bundle $P^r_{\cal F}M=\{j^r_0(f)|\, f:0\in U\subset
\R^n\to M, T_0f\in\op{Isom}(\R^n,T_{f(0)}M), Tf(T{\cal
F}_0)=T{\cal F}\},$ where ${\cal F}_0$ is the canonical regular
$p$-dimensional foliation of $\R^n$. The structure group of
$P^r_{\cal F}M$ is $G^r_{n,{\cal F}_0}=\{j^r_0(\zvf)|\, \zvf:0\in
U\subset \R^n\to\R^n, \zvf(0)=0, T_0\zvf\in\op{GL}(n,\R),
T\zvf(T{\cal F}_0)=T{\cal F}_0\}$, its action on $P^r_{\cal F}M$
is canonical. We call $P^r_{\cal F}M$ the principal bundle of {\sl
adapted $r$-frames} on $M$. For instance, $P^1_{{\cal
F}}M=:L_{\cal F}M$ is the bundle of adapted linear frames of $M$
with structure group \be G^1_{n,{\cal
F}_0}\simeq\op{GL}(n,q,\R)=\left\{\left(\begin{array}{cc}A&B\\0&D\end{array}\right):A\in\op{GL}(p,\R),B\in\op{gl}(p\times
q,\R),D\in\op{GL}(q,\R)\right\}.\ee Of course, the isotropy
subgroup of $[e_{n+1}]$ for the natural action of
\be\begin{array}{c}
\op{PGL}(n+1,q+1,\R)=\left\{\left(\begin{array}{ccc}A&B&h'\\0&D&h''\\0&\za''&a\end{array}\right):A\in\op{GL}(p,\R),B\in\op{gl}(p\times
q,\R),\right.\\\\\left.h'\in\R^p,D\in\op{GL}(q,\R),h''\in\R^q,\za''\in\R^{q*},a\in\R_0\right\}/\R_0\op{id}\end{array}\ee
on $\R P^n$ is \be\begin{array}{c}
H(n+1,q+1,\R)=\left\{\left(\begin{array}{ccc}A&B&0\\0&D&0\\0&\za''&a\end{array}\right):A\in\op{GL}(p,\R),B\in\op{gl}(p\times
q,\R),\right.\\\\\left.D\in\op{GL}(q,\R),\za''\in\R^{q*},a\in\R_0\right\}/\R_0\op{id}.\end{array}\ee

\begin{prop} Inclusion $I:H(n+1,\R)\to G^2_n$ of Proposition \ref{InclusionUsualContext} restricts to an inclusion
$I_{{\cal F}_0}:H(n+1,q+1,\R)\to G^2_{n,{\cal
F}_0}$.\label{InclusionAdapContext}\end{prop}

\begin{proof} It follows from the proof of Proposition \ref{InclusionUsualContext}
that the representative matrix of the tangent map at $0$ of the
smooth map $\zvf$ induced by an element of $H(n+1,q+1,\R)$ is
${\cal A}=\left(\begin{array}{cc}A & B\\ 0 & D\end{array}\right)$.
Hence the conclusion.
\end{proof}

We are now prepared to word the adapted version of Theorem
\ref{ProjStruReductions}.

\begin{theo}\label{ProjClassCartBundAdap} For any foliated manifold $(M,{\cal F})$, there exists a canonical
injection from the set of projective classes of adapted
connections $[\n_{{\cal F}}]$ into the set of reductions $P_{\cal
F}$ of the principal bundle $P^2_{\cal F}M$ of ${\cal F}$-adapted
second order frames on $M$ to structure group
$H(n+1,q+1,\R)\subset G^2_{n,{\cal F}_0}$.
\end{theo}

\begin{proof} The proof consists of three stages.\\

1. Let $\n_{\cal F}$ be an adapted connection of a foliated
manifold $(M,{\cal F})$. We will define the reduction $P_{\cal F}$
of $P^2_{\cal F}M$ to $H:=H(n+1,q+1,\R)$ by means of local
sections $\zs_{\za}$ of $P^2_{\cal F}M$ over open domains
$W_{\za}\subset M$ of adapted coordinates
$X_{\za}=(x_{\za},y_{\za})$ that form a cover $(W_{\za})$ of $M$.
Of course, the fiber $P_{{\cal F},m}$ of $P_{\cal F}$ at $m\in
W_{\za}$ is then $\zs_{\za}(m)\cdot H$, where $\cdot$ denotes the
action of $G^2_{n,{\cal F}_0}$ on $P^2_{\cal F}M$. The reduction
$P_{\cal F}=\cup_{m\in M}P_{{\cal F},m}$ is well-defined if and
only if the corresponding cocycle ${\frak
s}_{\za\zb}:W_{\za\zb}:=W_{\za}\cap W_{\zb}\to G^2_{n,{\cal
F}_0}$, which links the local sections,
$\zs_{\za}=\zs_{\zb}\cdot{\frak s}_{\zb\za}$, is valued in $H.$

If  $\zd^{i}_k$ is Kronecker's symbol and
$\zG^{i}_{\za;kl}\in\Ci(W_{\za})$ are Christoffel's symbols of
$\n_{\cal F}$, we set \be\zs_{\za}:W_{\za}\ni m\mapsto
(X^{i}_{\za}(m),\zd^{i}_k,-\zG^{i}_{\za;kl}(m))_{X_{\za}}\in
P^2_{{\cal F},m}M.\label{DefAdaCartBd}\ee Indeed, the image of $m$
is the package of partial derivatives that characterizes in the
coordinates $X_{\za}$ the $2$-jet at $0$ of the function $f:0\in
U\subset \R^n\to M$, $T_0f\in\op{Isom}(\R^n,T_{f(0)}M)$,
$Tf(T{\cal F}_0)=T{\cal F}$, with local form
$F^{i}_{\za}(Z)=X^{i}_{\za}(f(Z))=X^{i}_{\za}(m)+Z^{i}-\frac{1}{2}\zG^{i}_{\za;kl}(m)Z^kZ^l$,
where $Z\in U\subset\R^n$. In order to compare
$\zs_{\za}(m)=(X^{i}_{\za}(m),\zd^{i}_k,-\zG^{i}_{\za;kl}(m))_{X_{\za}}$
and
$\zs_{\zb}(m)=(X^{i}_{\zb}(m),\zd^{i}_k,-\zG^{i}_{\zb;kl}(m))_{X_{\zb}}$
for $m\in W_{\za\zb}$, we write $\zs_{\za}(m)$ using its
characterizing package of derivatives in the adapted coordinates
$X_{\zb}$. When applying formul$\ae\mbox{}$
(\ref{JetCharChangeCoord1}) and (\ref{JetCharChangeCoord2}), the
transformation law
$$\zG^{a}_{\za;kl}=\p_{X_{\za}^k}X_{\zb}^b\p_{X_{\za}^l}X_{\zb}^c\p_{X_{\zb}^d}X_{\za}^{a}\zG^d_{\zb;bc}+
\p^2_{X_{\za}^kX_{\za}^l}X_{\zb}^d\p_{X_{\zb}^d}X_{\za}^{a}$$ of
Christoffel's symbols, Equation (\ref{ActionChar}), as well as
Propositions \ref{InclusionUsualContext} and
\ref{InclusionAdapContext}, we get
$$\begin{array}{ll}\zs_{\za}(m)&=(X_{\zb}^{i}(m),\p_{X_{\za}^k}X_{\zb}^{i}(m),
-\zG^{i}_{\zb;bc}(m)\p_{X_{\za}^k}X_{\zb}^b(m)\p_{X_{\za}^l}X_{\zb}^c(m))_{X_{\zb}}\\
&=(X_{\zb}^{i}(m),\zd^{i}_k,-\zG^{i}_{\zb;kl}(m))_{X_{\zb}}\cdot(0,\p_{X_{\za}^k}X_{\zb}^{i}(m),0)\\
&=\zs_{\zb}(m)\cdot
\left[\left(\begin{array}{cc}\p_{X_{\za}}X_{\zb}(m)&0\\0&1\end{array}\right)\right]\\
&=\zs_{\zb}(m)\cdot {\frak s}_{\zb\za}(m).\end{array}$$ Since both
coordinate systems, $X_{\za}=(x_{\za},y_{\za})$ and
$X_{\zb}=(x_{\zb},y_{\zb})$, are adapted, we have
$y_{\zb}=y_{\zb}(y_{\za})$, so that ${\frak
s}_{\zb\za}:W_{\za\zb}\to H$.\\

2. We now prove that the just constructed reduction $P_{\cal F}$
of $P^2_{\cal F}M$ to $H$ does not depend on the considered
adapted connection $\n_{\cal F}$, but only on the projective class
$[\n_{\cal F}]$ of this connection. If $\n'_{\cal F}$ is a
projectively equivalent adapted connection, and if we set
$\p_{X_{\za}^{i}}=\p_{i}$, Equation (\ref{ProjEquiv}) entails that
$\n'_{{\cal F};\p_k}\p_{l}-\n_{{\cal
F};\p_k}\p_{l}=(\zG'^{i}_{\za;kl}-\zG_{\za;kl}^{i})\p_{i}=\za_k\p_l+\za_l\p_k.$
It follows that
$$\begin{array}{ll}\zs'_{\za}(m)&=(X^{i}_{\za}(m),\zd^{i}_k,-\zG'^{i}_{\za;kl}(m))_{X_{\za}}
=(X^{i}_{\za}(m),\zd^{i}_k,-\zG_{\za;kl}^{i}(m)-\zd^{i}_k\za_l(m)-\zd^{i}_l\za_k(m))_{X_{\za}}\\
&=(X^{i}_{\za}(m),\zd^{i}_k,-\zG^{i}_{\za;kl}(m))_{X_{\za}}\cdot(0,\zd^{i}_k,-\zd^{i}_k\za_l(m)-\zd^{i}_l\za_k(m))\\
&=\zs_{\za}(m)\cdot\left[\left(\begin{array}{cc}\zd^{i}_k&0\\\za(m)&1\end{array}\right)\right]=:\zs_{\za}(m)\cdot
h,\end{array}$$ where we used again Equation (\ref{ActionChar}).
As, in view of Proposition \ref{AdapProjEquiv},
$\za\in\zW^1(M,{\cal F})$ is foliated, we
have $\za_{\zi}=0$ and $h\in H$.\\

3. If the images $P_{\cal F}[\n_{\cal F}]$ and $P_{\cal
F}[\n'_{\cal F}]$ coincide, their fibers over any domain $W_{\za}$
of adapted coordinates $X_{\za}$ coincide. In particular, for any
$m\in W_{\za}$, there is a unique $h_{\za}(m)\in H$, such that
$\zs'_{\za}(m)=\zs_{\za}(m)\cdot h_{\za}(m).$ Hence,
$h_{\za}(m)=(0,{\cal A}^{i}_{\za;k}(m),-{\cal
A}^{i}_{\za;k}(m)\xi_{\za;l}(m)-{\cal
A}^{i}_{\za;l}(m)\xi_{\za;k}(m))$, see Proposition
\ref{InclusionUsualContext}, where $\xi_{\za;\zi}(m)=0$, see
Proposition \ref{InclusionAdapContext}. It easily follows that
${\cal A}^{i}_{\za;k}(m)=\zd^{i}_k$ and that
$\zG'^{i}_{\za;kl}(m)=\zG^{i}_{\za;kl}(m)+\zd^{i}_k\xi_{\za;l}(m)+\zd^{i}_l\xi_{\za;k}(m).$
Thus, for any $X,Y\in\op{Vect}(M)$, we have on $W_{\za}$,
$\n'_{{\cal F};X}Y-\n_{{\cal F};X}Y=\xi_{\za}(Y)X+\xi_{\za}(X)Y$.
Hence, $\xi_{\za}=\xi_{\zb}$ on $W_{\za\zb}$, and the
$\xi_{\za}\in\Ci(W_{\za},\R^{n*})$ define a unique differential
$1$-form $\xi\in\zW^1(M)$. Eventually, we get $[\n_{\cal
F}]=[\n'_{\cal F}]$.\end{proof}

\subsubsection{Foliated frame bundles}\label{FoliFramBund}

We next prove existence of a similar injection from projective
classes of foliated connections
into reductions of the ``foliated'' second order frame bundle.\\

Consider a foliated manifold $(M,{\cal F})$ and let ${\cal
U}=(U_i,f_i,g_{ij})$ be a H$\ae$fliger cocycle of ${\cal F}$ with
associated transverse manifold $N$. As $f_i$ is a submersion the
fibers (preimages) of which are parts of the leaves of ${\cal F}$,
the kernel of $T_{m}f_i:T_mU_i\to T_{f_i(m)}N_i$, $m\in U_i$, is
$\op{ker}T_mf_i=T_m{\cal F}$, and $N_mf_i:N_m(U_i,{\cal
F}):=T_mM/T_m{\cal F}\ni[v]\to (T_mf_i)(v)\in
T_{f_i(m)}N=N_{f_i(m)}(N_i,0)$ is a vector space isomorphism. Of
course, $(N_i,0)$ denotes the manifold $N_i$ endowed with its
canonical foliation by points. Actually, the {\sl normal functor}
$N$ is a functor between the category ${\cal FM}_q$ of codimension
$q$ foliated manifolds and smooth maps that preserve the
foliations, on one hand, and the category ${\cal FB}$ of foliated
fiber bundles, i.e. fiber bundles whose total space is foliated by
a foliation whose leaves are covering space of leaves on the base
space and bundle maps, on the other (see \cite{RW}). If confusion
with the transverse manifold $N$ is excluded, most authors denote
the {\sl normal bundle} $N(M,{\cal F})$ simply by $N$. Observe
also that $Nf_i$ is just the tangent map $Tf_i$ viewed as map
between normal bundles.

We now define the principal bundle $P^r(M,{\cal F})$, $r\in\N_0$,
of normal $r$th order frames associated with any object $(M,{\cal
F})\in \op{Obj}({\cal FM}_q)$. Remark first that each vector space
isomorphism $\R^q\to N_m=N_m(M,{\cal F})$, $m\in M$, implements a
normal linear frame $(n_1,\ldots,n_q)\in N^{\times q}_m$. In order
to obtain such isomorphisms, we consider the jets of {\sl
transverse smooth maps} $f: 0\in V\subset \R^q\to M$, such that
$\op{im}Tf\oplus T{\cal F}=TM$. Indeed, then $N_{z''}f:\R^q\ni
v\to [(T_{z''}f)(v)]\in N_{f(z'')}$, $z''\in V$, is a vector space
isomorphism. Hence, the set $P^r(M,{\cal F})$ (the better notation
$P^rN(M,{\cal F})$ is not prevailing) of normal $r$-frames is
defined by $P^r(M,{\cal F})=\{J^r_0(f)|\, f:0\in V\subset \R^q\to
M,\op{im}Tf\oplus T{\cal F}=TM\}$. The $r$-jet $J^r_0(f)$ at $0$
of a transverse function $f$ is the equivalence class of $f$ for
the following relation: two transverse functions $f$ and $g$ that
map a neighborhood $V\subset\R^q$ of $0$ into $M$ are equivalent
if and only if $f(0)=g(0)=:m$ and, for any submersion ${\frak
X}:m\in W\subset M\to \R^q$ that is constant along the leaves of
${\cal F}$, the components of the maps ${\frak F}:={\frak X}\circ
f$ and ${\frak G}:={\frak X}\circ g$ have the same partial
derivatives at $0$ up to order $r$. Of course, it suffices that
this condition be satisfied for one submersion. If $X=(x,y)$ is a
system of adapted coordinates of $M$ around $m$, we can choose
${\frak X}=y$. It is helpful to observe that (just as $Tf_i$, see
above) $T{\frak X}$ is a pointwise isomorphism of vector spaces
from $N$ onto $\R^q$, so that $T{\frak F}=T{\frak X}\circ Tf$. The
just defined space $P^r(M,{\cal F})$ is a principal bundle over
$M$ with structure group $G^r_q$ and projection
$\zp^r:J^r_0(f)\mapsto f(0)$. The right action is given by
$J^r_0(f)\cdot j^r_0(\zvf)=J^r_0(f\circ\zvf)$. Bundle $P^1(M,{\cal
F})=:L(M,{\cal F})$ for instance, is the principal bundle of
normal linear frames. Just as $N$ (see above), $P^r$ (or better
$P^rN$) is a functor between the categories ${\cal FM}_q$ and
${\cal FB}$. Let us mention that both functors are (prototypes of)
{\sl foliated natural functors} in the sense of \cite{RW}.

\begin{theo} For any foliated manifold $(M,{\cal F})$ of
codimension $q$, there exists an injection from the set of
projective classes of foliated connections $[\n(\cal F)]$ into the
set of reductions $P({\cal F})$ of $P^2(M,{\cal F})$ to structure
group $H(q+1,\R)\subset
G^2_q$.\label{ProjClassCartBundFol}\end{theo}

\begin{proof} The proof of this theorem is similar to that of
Theorem \ref{ProjClassCartBundAdap}. Hence, we put down only a
sketch of this proof.

If $\n({\cal F})$ is a foliated connection of a foliated manifold
$(M,{\cal F})$, the reduction $P({\cal F})$ of $P^2(M,{\cal F})$
to $H:=H(q+1,\R)$ is defined, over an open domain $W_{\za}\subset
M$ of adapted coordinates $X_{\za}=(x_{\za},y_{\za})$, by a local
section \be\zs_{\za}:W_{\za}\ni m\mapsto (y^{\frak
i}_{\za}(m),\zd^{\frak i}_{\frak k},-\zG^{\frak i}_{\za;{\frak
k}{\frak l}}(m))_{X_{\za}}\in P^2_m(M,{\cal
F}),\label{DefFolCarBd}\ee where the $\zG^{\frak i}_{\za;{\frak
k}{\frak l}}\in\Ci(W_{\za})$ are Christoffel's symbols of
$\n({\cal F})$. A similar argument than in the adapted case, again
allows checking that the cocycle ${\frak s}_{\zb\za}$, which links
$\zs_{\za}$ and $\zs_{\zb}$, is valued in subgroup $H$. Also
invariance of the reduction for a change of foliated connection
within the same projective class, as well as injectivity of the
just defined mapping between projective classes and reductions,
can be verified as above.\end{proof}

\subsubsection{Projections}

It is a well-known fact (see above, adapted and foliated vector
fields, adapted and foliated connections) that adapted objects
induce (usually) foliated objects. In this subsection, we describe
canonical projections from an adapted frame bundle $P^r_{\cal F}M$
(resp. {\sl adapted Cartan bundle} $P_{\cal F}$) onto the
corresponding foliated frame bundle $P^r(M,{\cal F})$ (resp. {\sl
foliated Cartan bundle} $P({\cal F})$).\\

We denote by $p_{\cal F}^{r}$ (resp. $p^{r}({\cal F})$), $r\ge 1$,
the canonical projection $p_{\cal F}^{r}:P_{\cal F}^rM\ni
j^r_0(f)\mapsto j^{r-1}_0(f)\in P_{{\cal F}}^{r-1}M$ (resp.
$p^{r}({\cal F}):P^r(M,{\cal F})\ni J^r_0(f)\mapsto
J^{r-1}_0(f)\in P^{r-1}(M,{\cal F})$). Furthermore, if $f:0\in
U\subset \R^n\to M$, $T_0f\in\op{Isom}(\R^n,T_{f(0)}M)$,
$Tf(T{\cal F}_0)=T{\cal F}$, and if $i_q:\R^q\ni z''\mapsto
(0,z'')\in \R^n$, then, obviously, $f\circ i_q:0\in V\subset
\R^q\to M$, $\op{im}T(f\circ i_q)\oplus T{\cal F}=TM.$ Since two
foliation preserving locally defined diffeomorphisms $f$ and $g$
that have the same jets $j^r_0(f)=j^r_0(g)$, induce two transverse
maps $f\circ i_q$ and $g\circ i_q$, such that $J^r_0(f\circ
i_q)=J^r_0(g\circ i_q)$, there is a canonical projection $^{\cal
F}\!p^{r}:P^r_{\cal F}M\ni j^r_0(f)\mapsto J^r_0(f\circ i_q)\in
P^r(M,{\cal F})$. Observe that if
$Z=(z',z'')\in\R^p\times\R^q=\R^n$ (resp. $X=(x,y)$) are adapted
coordinates in $(\R^n,{\cal F}_0)$ (resp. $(M,{\cal F})$ around
$f(0)$), jet $j^r_0(f)$ is characterized by the derivatives
$$\p_{Z}^{\za}(X^{i}(f(Z)))(0),\vert\za\vert\le r,i\in\{1,\ldots,n\},$$ whereas jet $J^r_0(f\circ i_q)$ is
characterized by $$\p_{z''}^{\za}(y^{\frak
i}(f(0,z'')))(0)=\p_{z''}^{\za}(y^{\frak
i}(f(Z)))(0,0),\vert\za\vert\le r,{\frak i}\in\{1,\ldots,q\}.$$

\begin{prop} For any foliated manifold $(M,{\cal F})$ endowed with
an adapted projective structure and the induced foliated
projective structure, projection $^{\cal F}\!p^2:P^2_{\cal F}M\to
P^2(M,{\cal F})$ restricts to a projection $^{\cal F}\!{\frak
p}^2:P_{\cal F}\to P({\cal F})$, and the diagram
$$\begin{array}{ccccc}&P_{\cal F}\subset P^2_{\cal F}M&\stackrel{p^2_{\cal F}}{\to}& L_{\cal F}M&\\\\
^{\cal F}\!{\frak p}^2 &\downarrow &&\downarrow & ^{\cal
F}\!p^1\\\\&P({\cal F})\subset P^2(M,{\cal F})&\stackrel{p^2({\cal
F})}{\to}& L(M,{\cal F})\end{array}$$ is
commutative.\label{ProjCartBundAdapCartBundFoli}\end{prop}

\begin{proof} It suffices to prove that $^{\cal F}\!{\frak p}^2$ maps $P_{\cal F}$ into $P({\cal
F})$. Consider a point of $P_{{\cal F},m}$, $m\in M$, i.e., in
adapted coordinates $X=(x,y)$ around $m$, see Equation
(\ref{DefAdaCartBd}), a point
$$(X^{i}(m),\zd^{i}_k,-\zG^{i}_{{\cal F},kl}(m))_{X}\cdot
(0^{i},{\cal A}^{i}_k,-{\cal A}^{i}_k\za_l-{\cal A}^{i}_l\za_k),$$
where the element of $G^2_{n,{\cal F}_0}$ is induced by a member
$$\left(\begin{array}{cc}{\cal A}&0\\\za & 1\end{array}\right),
{\cal A}=\left(\begin{array}{cc}A&B\\0&D\end{array}\right),
\za=(0,\za''),$$ of $H(n+1,q+1,\R).$ When using Equation
(\ref{ActionChar}), the above description of $^{\cal F}\!p^r$ in
terms of packages of derivatives, the local characterization of an
adapted connection, see Equation (\ref{AdapConnLocal1}), as well
as Proposition \ref{AdapConnIndFolConn}, we see that the
considered point of $P_{{\cal F},m}$ is mapped by $^{\cal
F}\!{\frak p}^2$ to
$$(y^{\frak i}(m),D^{\frak i}_{\frak k},-D^{\frak i}_{\frak
k}\za''_{\frak l}-D^{\frak i}_{\frak l}\za''_{\frak k}-\zG({\cal
F})^{\frak i}_{{\frak a}{\frak b}}(m)D^{\frak a}_{\frak k}D^{\frak
b}_{\frak l})_X=(y^{\frak i}(m),\zd^{\frak i}_{\frak k},-\zG({\cal
F})^{\frak i}_{{\frak k}{\frak l}}(m))_X\cdot(0^{\frak i},D^{\frak
i}_{\frak k},-D^{\frak i}_{\frak k}\za''_{\frak l}-D^{\frak
i}_{\frak l}\za''_{\frak k}).$$ It then follows directly from
Equation (\ref{DefFolCarBd}) and Proposition
\ref{InclusionUsualContext} that $^{\cal F}\!{\frak p}^2$ is
valued in $P(\cal F)$.
\end{proof}

\section{Lift of adapted and foliated symbols}

Below, we study adapted and foliated symbols, as well as their
lifts to the Cartan fiber bundles $P_{\cal F}$ and $P(\cal F).$
Investigations are again similar in both settings. Whereas we
detailed above the adapted situation, we describe below especially
the foliated case.

\subsection{Foliated differential operators and symbols}

We have already mentioned, see Subsection \ref{FoliFramBund}, that
$N$, $LN$, and more generally $P^rN$, $r\in\N_0$, are (covariant)
foliated natural functors, i.e. (regular) functors $F:{\cal
FM}_q\to {\cal FB}$, such that for any morphism $f:(M_1,{\cal
F}_1)\to (M_2,{\cal F}_2)$, morphism $F(f):F(M_1,{\cal F}_1)\to
F(M_2,{\cal F}_2)$ covers $f$ and is pointwise a diffeomorphism.

For instance, if $(M,{\cal F})$ is defined by means of a cocycle
${\cal U}=(U_i,f_i,g_{ij})$, submersion $f_i:(U_i,{\cal F})\to
(N_i,0)$ (resp. diffeomorphism $g_{ji}:(N_{ij},0) \to (N_{ji},0)$)
is a morphism of ${\cal FM}_q$, and the corresponding morphism
$N(f_i)$ (resp. $N(g_{ji})$) is a pointwise isomorphism \be
N_m(f_i): N_m(U_i,{\cal F})\to N_{f_i(m)}(N_i,0)=
T_{f_i(m)}N_i,\label{NormBundPullBack}\ee $m\in U_i$ (resp. \be
N_m(g_{ji})=T_m(g_{ji}):N_m(N_{ij},0)= T_mN_{ij}\to
N_{g_{ji}(m)}(N_{ji},0)= T_{g_{ji}(m)}N_{ji},\label{GluiF_N}\ee
$m\in N_{ij}$). It is easily checked that $(N(U_i,{\cal
F}),N(f_i), N(g_{ij}))$, and more generally $ ( F(U_i,{\cal F}),
F(f_i),\linebreak F(g_{ij} ) ), $ is a cocycle that defines a
foliation ${\cal F}_F$ on the total space $F(M,{\cal F})$ (and
that ${\cal F}_F$ is independent of ${\cal U}$). Hence, the name
``foliated natural bundle''. Further, it follows from Equation
(\ref{NormBundPullBack}) that \be N(U_i,{\cal F})\simeq
f_i^*N(N_i,0)=f_i^* T{N_i}.\label{NormBundPullBack2}\ee

More generally, $F_{\cal U} = \amalg_i f_{i}^{*}F(N_{i},0)/\sim,$
where $(i,m',v') \sim (j,m'',v'')$ if and only if $m' = m''$ and
$v'' = F(g_{ji})(v')$, is a well-defined fibre bundle over $M$.
The projections \be f_{i}^{*}F(N_{i},0) \to
F(N_{i},0)\label{InduFoli}\ee define a foliation on $F_{\cal U}$.
It is obvious, see preceding equations, that (the foliated) bundle
$F(M,\cal F)$ is isomorphic to (the foliated) bundle $F_{\cal U}$
and that this isomorphism is foliation preserving.

The ``mental picture'' of foliation ${\cal F}_F$ induced on
$F(M,{\cal F})$ is clear from Equation (\ref{InduFoli}). In
particular, foliation ${\cal F}_F$ has the same dimension as
foliation ${\cal F}$ and its leaves project onto the leaves of
${\cal F}$.

Let us also recall (see \cite{RW}) that a {\sl foliated geometric
structure} is a foliated subbundle of a foliated natural bundle
$F(M,{\cal F})$, i.e. a subbundle (in particular a section) the
total space of which is saturated for foliation ${\cal F}_F$ (``it
contains as many leaves as can reasonably be expected'').

\begin{defi} A {\sl foliated differential operator} of a foliated manifold $(M,{\cal F})$
$($where ${\cal F}$ is of dimension $p$ and codimension $q$$)$ is
an endomorphism $D\in\op{End}_{\R}(\Ci(M,{\cal F}))$ that reads in
any system of adapted coordinates
$(x,y)=(x^1,\ldots,x^p,y^1,\ldots,y^q)$ over any open subset
$U\subset M$,
$$D\vert_U = \sum_{\vert \zg \vert \le k} D_{\zg}\,[\p_{y^{1}}]^{\zg^1}\ldots[\p_{y^q}]^{\zg^q},$$
where $k\in\N$ is independent of the considered adapted chart and
where the coefficients $D_{\zg}\in\Ci(U,{\cal F})$ are locally
defined foliated functions. The smallest possible integer $k$ is
called the {\sl order} of operator $D$.
\end{defi}

We denote by ${\cal D}(M,{\cal F})$ (resp. ${\cal D}^k(M,{\cal
F})$) the space of all foliated differential operators (resp. all
foliated differential operators of order $\le k$). Of course, the
usual filtration \be {\cal D}(M,{\cal F})= \cup_{k \in {\N}}{\cal
D}^k(M,{\cal F})\label{FiltFoliDiffOper}\ee holds true.

\begin{defi} The graded space ${\cal S}(M,{\cal F})$ associated with the filtered space ${\cal D}(M,{\cal F})$,
$${\cal S}(M,{\cal F})= \oplus_{k \in \N} {\cal S}^k (M, {\cal F})=\oplus_{k \in \N}{\cal D}^k (M,{\cal F})/{\cal D}^{k-1}
(M, {\cal F}),$$ is the space of {\sl foliated symbols}.
\end{defi}

It is easily checked that the well-known vector space isomorphism
between the spaces of symbols of degree $k$ and of symmetric
contravariant $k$-tensor fields, extends to the foliated setting,
$${\cal S}^k (M, {\cal F})\simeq \Gamma(S^kN(M,{\cal F});{\cal F}_{S^kN}),$$
where the {\small RHS} denotes the space of foliated sections of
the foliated natural bundle $S^kN(M,{\cal F})$ (see above,
foliated geometric structures). Below, we identify these two
spaces.

\begin{theo}\label{LiftFoliSymbTheo} Let $(M,{\cal F})$ be a foliated manifold of codimension $q$
endowed with a foliated projective structure $[\n({\cal F})]$, and
denote by $P({\cal F})$ the corresponding reduction of
$P^2(M,{\cal F})$ to $H(q+1,\R)\subset G^2_q$. The following
canonical vector space isomorphisms hold: \be \sim : {\cal
S}^k(M,{\cal F})= \Gamma(S^kN(M,{\cal F});{\cal F}_{S^kN})\ni
s\mapsto \tilde{s}\in\Ci(LN(M,{\cal F}),S^k\R^q;{\cal
F}_{LN})_{\op{GL}(q,\R)}\label{LiftFoliSymb1}\ee \be \wedge:{\cal
S}^k(M,{\cal F})= \Gamma(S^kN(M,{\cal F});{\cal F}_{S^kN})\ni
s\mapsto \hat{s}=\tilde{s}\circ p^2({\cal F})\in\Ci(P({\cal
F}),S^k\R^q;{\cal
F}_{P^2N})_{H(q+1,\R)}.\label{LiftFoliSymb2}\ee\end{theo}

\begin{proof} 1. Observe first that the foliated natural vector bundle $S^kN(M,{\cal F})$
is associated with the foliated natural principal bundle
$LN(M,{\cal F})$ of normal linear frames: $S^kN(M,{\cal
F})=LN(M,{\cal F})\times_{\op{GL}(q,\R)} S^k\R^q$. Hence, only the
foliated aspect of Isomorphism (\ref{LiftFoliSymb1}) has to be
explained. Consider a section $s\in\Gamma(S^kN(M,{\cal F});{\cal
F}_{S^kN})$ and a normal linear frame $u_m=([v_1],\ldots,[v_q])\in
LN_m(M,{\cal F})$. Isomorphism $\sim$ is of course defined by
$\tilde{s}(u_m)=(s^{i_1\ldots i_k}(m))\in S^k\R^q$, where the
{\small RHS} is made up by the components of
$s_m=\sum_{i_1\le\ldots\le i_k}s^{i_1\ldots i_k}(m)
[v_{i_1}]\vee\ldots\vee[v_{i_k}]$ in the induced linear frame of
$S^kN_m(M,{\cal F})$. Hence, the $\op{GL}(q,\R)$-equivariance of
$\tilde{s}$ is obvious. But this function is also foliated, i.e.
locally constant along the leaves of ${\cal F}_{LN}$. Indeed, let
$(U_i,f_i,g_{ij})$ be a defining cocycle of ${\cal F}$, and let
$u'_{m'}\in LN_{m'}(M,{\cal F})$ be a normal linear frame on the
same local leave of ${\cal F}_{LN}$ than $u_m$; the leaves of
${\cal F}_{LN}$ are locally defined by the projections \be
f_i^*LN(N_i,0)=f_i^*LTN_i\to LTN_i.\label{FoliationLN}\ee Since
section $s$ is foliated and as the local leaves of the
corresponding foliation ${\cal F}_{S^kN}$ are defined by the
projections $f_i^*S^kN(N_i,0)=f_i^*S^kTN_i\to S^kTN_i$, it is
clear that the tensors $s_m$ and $s_{m'}$ have the same components
in the frames $u_m$ and $u'_{m'}$ respectively, so that
$\tilde{s}(u_m)=\tilde{s}(u'_{m'})$. A similar argument shows that
to any foliated equivariant function is associated a foliated section. \\

2. We will show in Point 3 that the spaces of foliated equivariant
functions on $LN(M,{\cal F})$ and on $P({\cal F})$, see Equations
(\ref{LiftFoliSymb1}) and (\ref{LiftFoliSymb2}), are isomorphic.
This is a foliated variant of a result that has already been
proven in \cite{MR}.

Let us recall that the action $\tilde{\zr}$ of $H(q+1,\R)\simeq
\op{GL}(q,\R)\rtimes \R^{q*}$ on $S^k\R^q$ is induced by the
action $\zr$ of $\op{GL}(q,\R)$ on $S^k\R^q$:
\be\tilde{\zr}\left[\left(\begin{array}{cc}{\cal A}& 0\\
\za & a\end{array}\right)\right]=\zr(\frac{\cal
A}{a}).\label{ActiH}\ee

In order to understand that the target space of Equation
(\ref{LiftFoliSymb2}) makes sense, observe that $P({\cal F})$ is a
foliated subbundle of the foliated natural bundle $P^2N(M,{\cal
F})$. Indeed, first it is clear that if $\zf:(M_1,{\cal F}_1)\to
(M_2,{\cal F}_2)$ is a morphism of category ${\cal FM}_q$, then
the corresponding morphism $P^rN(\zf)$ of category ${\cal FB}$ is
defined by $P^rN(\zf):P^rN(M_1,{\cal F}_1)\ni J^r_0(f)\mapsto
J^r_0(\zf\circ f)\in P^rN(M_2,{\cal F}_2)$; in particular,
$P^rN(f_i):P^rN(U_i,{\cal F})\ni J^r_0(f)\mapsto j^r_0(f_i\circ
f)\in P^rTN_i$, with self-explaining notations. Furthermore,
since, in adapted coordinates $X=(x,y)$ over an open subset
$W\subset M$, we have $\zs(x,y)=(y^{\frak i},\zd^{\frak i}_{\frak
k},-\zG({\cal F})^{\frak i}_{{\frak k}{\frak l}}(y))\in
P_{(x,y)}({\cal F}),$ see Equation (\ref{DefFolCarBd}), the local
section $\zs$ of $P({\cal F})$ is constant along any local leaf of
${\cal F}$ in $W$. Hence, $\zs$ is valued in a leaf of ${\cal
F}_{P^2N}$. Eventually, the action by an element $h=j^2_0(\zvf)\in
H(q+1,\R)\subset G^2_q$ maps a local leaf of ${\cal F}_{P^2N}$
into another local leaf. As a matter of fact, if $u^2=J^2_0(f)$
and $u'^2=J^2_0(f')$ belong to the same local leaf of ${\cal
F}_{P^2N}$, we have $j^2_0(f_i\circ
f)=P^2N(f_i)(u^2)=P^2N(f_i)(u'^2)=j^2_0(f_i\circ f')$. But then,
$P^2N(f_i)(u^2\cdot h)=j^2_0(f_i\circ f\circ \zvf)=j^2_0(f_i\circ
f'\circ \zvf)=P^2N(f_i)(u'^2\cdot h)$. Thus, $P({\cal F})$
is actually a foliated subbundle.\\

3. Observe first that mapping $p:=p^2({\cal F}):J^2_0(f)\in
P({\cal F})\subset P^2N(M,{\cal F})\to J^1_0(f)\in LN(M,{\cal F})$
is surjective. Indeed, the fiber $LN_m(M,{\cal F})$, $m\in M$, is
equivalent to $\op{GL}(q,\R)$. On the other hand, in adapted
coordinates $X=(x,y)$ around $m$, the projection of the
corresponding fiber $P_m({\cal F})$ of $P({\cal F})$ is made up,
see Equation (\ref{DefFolCarBd}), by the elements \be
p(\zs(m)\cdot h) = p((y^{\frak i}(m),\zd^{\frak i}_{\frak
k},-\zG({\cal F})^{\frak i}_{{\frak k}{\frak l}}(m))_{X}\cdot
(0^{\frak i},{\cal A}^{\frak i}_{\frak k},-{\cal A}^{\frak
i}_{\frak k}\za_{\frak l}-{\cal A}^{\frak i}_{\frak l}\za_{\frak
k}))=(y^{\frak i}(m),{\cal A}^{\frak i}_{\frak
k})_X,\label{ProjPFeui}\ee where $h$ runs through $H(q+1,\R)$, so
that ${\cal A}$ runs through $\op{GL}(q,\R)$.

For any ${\frak f}\in\Ci(LN(M,{\cal F}),S^k\R^q;{\cal
F}_{LN})_{\op{GL}(q,\R)}$, we now set $\hat{\frak f}={\frak
f}\circ p\in\Ci(P({\cal F}),S^k\R^q).$ This map $\hat{\frak f}$ is
$H(q+1,\R)$-equivariant. Actually, for any $u_m^2\in P_m({\cal
F})$, $u_m^2=(y^{\frak i}(m),{\cal B}^{\frak i}_{\frak k},{\cal
T}^{\frak i}_{{\frak k}{\frak l}})_X$, and any $h=(0^{\frak
i},{\cal A}^{\frak i}_{\frak k},-{\cal A}^{\frak i}_{\frak
k}\za_{\frak l}-{\cal A}^{\frak i}_{\frak l}\za_{\frak k})\in
H(q+1,\R)$, we have \be\begin{array}{c} f(p(u_m^2\cdot
h))=f((y^{\frak i}(m),{\cal B}^{\frak i}_{\frak a}{\cal A}^{\frak
a}_{\frak k})_X)=f((y^{\frak i}(m),{\cal B}^{\frak i}_{\frak
k})_X\cdot(0^{\frak i},{\cal A}^{\frak i}_{\frak k}))\\=\zr({\cal
A}^{-1})(f(p(u_m^2)))=\tilde{\zr}(h^{-1})(f(p(u_m^2))),\end{array}\label{EquiFuncPFeui}\ee
in view of Equation (\ref{ActiH}). Eventually, $\hat{\frak f}\in
\Ci(P({\cal F}),S^k\R^q;{\cal F}_{P^2N})_{H(q+1,\R)}.$ Indeed, if
$u^2=J^2_0(f)$ and $u'^{2}=J^2_0(f')$ are two points on the same
local leaf of ${\cal F}_{P^2N}$ in $P({\cal F})$, see Point 2,
then $f(p(u^2))=f(J^1_0(f))=f(J^1_0(f'))=f(p(u'^2))$, since $f$ is
locally constant along the leaves of ${\cal F}_{LN}$.

It is clear that $\wedge:{\frak f}\mapsto \hat{\frak f}$ is linear
and injective (since $p$ is surjective). Map $\wedge$ is also
surjective. Indeed, any function $g\in \Ci(P({\cal
F}),S^k\R^q;{\cal F}_{P^2N})_{H(q+1,\R)}$ factors through
$LN(M,{\cal F})$, i.e. $g={\frak g}\circ p$. Note first that it
follows from Equations (\ref{ProjPFeui}) and (\ref{EquiFuncPFeui})
that for any $u^1\in LN(M,{\cal F})$ and any ${\cal A}\in
\op{GL}(q,\R)$, there is $u^2\in P({\cal F})$ and $h\in
H(q+1,\R)$, such that $p(u^2)=u^1$, ${\cal A}$ is the upper left
submatrix of $h$, and $p(u^2\cdot h)=u^1\cdot {\cal A}$. Hence,
${\frak g}(u^1\cdot {\cal A})=g(u^2\cdot
h)=\tilde{\zr}(h^{-1})g(u^2)=\zr({\cal A}^{-1}){\frak g}(u^1)$,
and ${\frak g}$ is $\op{GL}(q,\R)$-equivariant. It is also
well-defined, since, if $u^2,u'^2\in P({\cal F})$ project both
onto $u^1$, we have $u'^2=u^2\cdot h$, $h\in H(q+1,\R)$, and
$u^1=u^1\cdot {\cal A}$. Thus, ${\cal A}=\op{id}$, $g(u^2\cdot
h)=g(u^2)$, and ${\frak g}\in\Ci(LN(M,{\cal
F}),S^k\R^q)_{\op{GL}(q,\R)}$. In order to prove that ${\frak g}$
is foliated for ${\cal F}_{LN}$, observe that, as the leaves of
${\cal F}_{P^rN}$, $r\in\N_0$, are locally defined as the fibers
of submersion \be P^rN(f_i):P^rN(U_i,{\cal F})\ni J^r_0(f)\mapsto
j^r_0(f_i\circ f)\in P^rTN_i,\label{FoliationPrN}\ee the local
leaves are made up by the jets $J^r_0(f)\neq J^r_0(f')$ of those
transverse functions $f$ and $f'$, the last $q$ adapted
coordinates of which have the same derivatives at $0$ up to order
$r$, but that map $0$ to $m:=f(0)\neq f'(0)=:m'$ on the same local
leaf of ${\cal F}$ in $U_i$, see definitions of the jets $J^r_0$
and $j^r_0$. Equation (\ref{ProjPFeui}) then entails that two jets
$u^1=J^1_0(f)$ and $u'^1=J^1_0(f')$ on the same local leaf of
${\cal F}_{LN}$ are the projections of two jets $u^2$ and $u'^2$
of $P({\cal F})$ on the same local leaf of ${\cal F}_{P^2N}$.
Hence, ${\frak g}$ is foliated.\end{proof}

It is interesting to observe that (the mental picture associated
with) foliation ${\cal F}_{LN}$ is of course the same,
irrespective of the fact it is defined by Equation
(\ref{FoliationLN}) or by Equation (\ref{FoliationPrN}).

\subsection{Adapted differential operators and symbols}

As aforementioned, in order to limit the length of this paper, we
confine ourselves in the adapted case to a description of the main
points. Hence, we refrain for instance to give a general
description of {\sl adapted natural functors}, see \cite{RW}, but
provide examples of such functors.

Let $(M,{\cal F})$ be a foliated manifold of dimension $p$ and
codimension $q$ and denote by $(U_i,f_i,g_{ij})$ a cocycle of
${\cal F}$. Then, for any $r\in\N_0$,
$\underline{P}^rU_i:=P^r_{\cal F}U_i$, \be
\underline{P}^r(f_i):P^r_{\cal F}U_i\ni j^r_0(f)\mapsto
j^r_0(f_i\circ f\circ i_q)\in
P^r_0N_i=P^rN_i,\label{FoliationPrBar}\ee and
$\underline{P}^r(g_{ij}):P^rN_{ji}\ni j^r_0(g)\mapsto
j^r_0(g_{ij}\circ g)\in P^rN_{ij}$, form a cocycle that defines a
foliation ${\cal F}_{\underline{P}^r}$ on $P^r_{\cal F}M$.
Consider now an adapted atlas of $(M,{\cal F})$ and take an
adapted chart $(U,\zf)$, $\zf=(\zf_1,\zf_2):U\to\R^p\times\R^q.$
The map \be \zf_{2*}:L_{\cal
F}U\ni(v_1,\ldots,v_p,v_{p+1},\ldots,v_{p+q})\to
(T\zf_2(v_{p+1}),\ldots, T\zf_2(v_{p+q}))\in
L\R^q\label{FoliationLBar}\ee is a submersion that defines a
foliation ${\cal F}_{\underline{L}}$ on $L_{\cal F}M$. It is clear
that foliation ${\cal F}_{\underline{P}^1}$, defined by Equation
(\ref{FoliationPrBar}), and foliation ${\cal F}_{\underline{L}}$,
defined by Equation (\ref{FoliationLBar}), coincide. Moreover, the
leaves of ${\cal F}_{\underline{P}^1}={\cal F}_{\underline{L}}$
project onto the leaves of ${\cal F}$ (since
$\underline{P}^1(f_i)$ and $\zf_{2*}$ are bundle maps over $f_i$
and $\zf_2$, respectively) and the dimension of foliation ${\cal
F}_{\underline{L}}$ is $p+np$ (no restrictions imposed, neither on
$v_1,\ldots,v_p\in T{\cal F}$, nor on the tangential parts of
$v_{p+1},\ldots,v_{p+q}$). The ``mental picture'' of ${\cal
F}_{\underline{L}}$ follows.

Foliation ${\cal F}$ similarly induces a foliation ${\cal F}_{T}$
on $TM$ (defined for instance by means of $T\zf_2:TU\to\R^q$).
Observe that adapted vector fields of $M$, see Equation
(\ref{AdapVectFiel}), coincide with sections of $TM$ that are
foliated for ${\cal F}_T$: $\op{Vect}_{\cal F}(M)=\zG(TM;{\cal
F}_T).$

Let us also mention that {\sl adapted functions} are just foliated
functions: $\Ci_{\cal F}(M):=\Ci(M,{\cal F}).$

\begin{defi} An {\sl adapted differential operator} of a foliated manifold $(M,{\cal F})$
$($where ${\cal F}$ is of dimension $p$ and codimension $q$$)$ is
an endomorphism $D\in\op{End}_{\R}(\Ci(M))\cap
\op{End}_{\R}(\Ci_{\cal F}(M))$ that reads in any system of
adapted coordinates $(x,y)=(x^1,\ldots,x^p,y^1,\ldots,y^q)$ over
any open subset $U\subset M$,
$$D\vert_U = \sum_{\vert \zg \vert \le k} D_{\zg}\,\p_{x^{1}}^{\zg^1}\ldots\p_{x^p}^{\zg^p}\p_{y^1}^{\zg^{p+1}}\ldots\p_{y^q}^{\zg^{p+q}},$$
where $k\in\N$ is independent of the considered adapted chart,
where $D_{\zg}\in\Ci(U)$, and where the coefficients $D_{\zg}$
with $\zg^1=\ldots=\zg^p=0$ are locally defined adapted functions.
The smallest possible integer $k$ is called the {\sl order} of
operator $D$.
\end{defi}

We denote by ${\cal D}_{\cal F}(M)$ the filtered space of all
adapted differential operators on $(M,{\cal F})$. The
corresponding graded space ${\cal S}_{\cal F}(M)$ is the space of
{\sl adapted symbols} on $(M,{\cal F})$. Of course, ${\cal
S}_{\cal F}^k (M)\simeq \Gamma_{\cal F}(S^kTM),$ where the {\small
RHS} denotes the space of adapted sections of $S^kTM$. The
definition of adapted sections of $S^kTM$ is clear in view of the
definitions of adapted vector fields and adapted differential
operators. Furthermore, as in the case of the tangent bundle,
foliation ${\cal F}$ induces a foliation ${\cal F}_{S^kT}$ on
$S^kTM$ (defined by means of the extension
$(T\zf_2)^{\otimes}:S^kTU\to S^k\R^q$ of $T\zf_2=TU\to\R^q$) and
adapted sections of $S^kTM$ coincide with sections of this bundle
that are foliated for ${\cal F}_{S^kT}$: $S^k_{\cal
F}(M)\simeq\zG_{\cal F}(S^kTM)=\zG(S^kTM;{\cal F}_{S^kT}).$ Below,
we identify the first two of the preceding spaces.

\begin{theo}\label{LiftAdapSymbTheo} Let $(M,{\cal F})$ be a foliated manifold of codimension $q$
endowed with an adapted projective structure $[\n_{\cal F}]$, and
denote by $P_{\cal F}$ the corresponding reduction of $P^2_{\cal
F}M$ to $H(n+1,q+1,\R)\subset G^2_{n;{\cal F}_0}$. We then have
the following canonical vector space isomorphisms: \be \sim :{\cal
S}^k_{\cal F}(M)=\Gamma_{\cal F}(S^kTM)\ni s\mapsto
\tilde{s}\in\Ci_{\cal F}(L_{\cal
F}M,S^k\R^n)_{\op{GL}(n,q,\R)}\label{LiftAdapSymb1}\ee \be
\wedge:{\cal S}^k_{\cal F}(M)=\Gamma_{\cal F}(S^kTM)\ni s\mapsto
\hat{s}=\tilde{s}\circ p^2_{\cal F}\in\Ci_{\cal F}(P_{\cal
F},S^k\R^n)_{H(n+1,q+1,\R)}.\label{LiftAdapSymb2}\ee\end{theo}

The proof of this theorem is on the same lines than that of
Theorem \ref{LiftFoliSymbTheo}. Let us explain the meaning of
$\Ci_{\cal F}$ in Equations (\ref{LiftAdapSymb1}) and
(\ref{LiftAdapSymb2}). In the following, we denote by $p_{n,q}$
canonical projections, such as
$p_{n,q}:\op{GL}(n,q,\R)\to\op{GL}(q,\R)$, $p_{n,q}:S^k\R^n\to
S^k\R^q$, ... Adapted functions $f\in\Ci_{\cal F}(P^r_{\cal
F}M,S^k\R^n)$ are then the functions $f\in\Ci(P^r_{\cal
F}M,S^k\R^n)$, for which $p_{n,q}\circ f\in\Ci(P^r_{\cal
F}M,S^k\R^q;{\cal F}_{\underline{P}^r})$ is foliated for ${\cal
F}_{\underline P^r}$. If we set $k=0$, we get that adapted
functions coincide with foliated functions, see above.

\subsection{Projections}

Adapted symbols project onto foliated symbols. Indeed, we have the

\begin{prop} For any foliated manifold $(M,{\cal F})$, there is a canonical degree-preserving
projection $^{\cal F}\!\zp: {\cal S}^k_{\cal F}(M)\to {\cal
S}^k(M,{\cal F})$, $k\in\N$. If the considered foliated manifold
is endowed with an adapted and the corresponding foliated
projective structures, and if $\,^{\cal F}\!\tilde{\zp}$ $($resp.
$^{\cal F}\!\hat{\zp}$$)$ denotes projection $^{\cal F}\!\zp$ read
through the isomorphisms $\sim$ $($resp. $\wedge$$)$ detailed in
Theorems \ref{LiftFoliSymbTheo} and \ref{LiftAdapSymbTheo}, we
have, for any symbol $s\in {\cal S}^k_{\cal F}(M)$, $k\in\N$, \be
\left(^{\cal F}\!\tilde{\zp}\tilde{s}\right)\circ ^{\cal
F}\!\!p^1=p_{n,q}\circ \tilde{s}\quad (\mbox{resp. }\left(^{\cal
F}\!\hat{\zp}\hat{s}\right)\circ ^{\cal F}\!\!{\frak
p}^2=p_{n,q}\circ \hat{s}).\label{PiFoliTildeHat}\ee
\end{prop}
\medskip

\begin{proof} As usual, we denote by $n$ the dimension of $M$ and by $p$ (resp. $q$) the dimension (resp.
codimension) of ${\cal F}$. Since $^{\cal F}\!p^r:P^r_{\cal F}M\ni
j^r_0(f)\mapsto J^r_0(f\circ i_q)\in P^r(M,{\cal F})$, as
$j^1_0(f)$, $f(0)=m$, corresponds to the basis
$(v_1,\ldots,v_n)\in (T_mM)^{\times n}$,
$$v_i=\sum_j\p_{Z^{i}}(X^j(f(Z)))(0)\p_{X^j},$$ where $Z=(z',z'')$ are canonical coordinates in $\R^n=\R^p\times\R^q$ and $X=(x,y)$ are
adapted coordinates in $M$ around $m$, and as $J^1_0(f\circ i_q)$
corresponds to the basis $(n_1,\ldots,n_q)\in (N_m)^{\times q},$
$$n_{\frak i}=\sum_{\frak j}\p_{z''^{\frak i}}(y^{\frak
j}(f(Z)))(0)[\p_{y^{\frak j}}],$$ we see that $^{\cal
F}\!p^1(v_1,\ldots,v_n)=([v_{p+1}],\ldots,[v_{p+q}])$.

Consider now $s\in {\cal S}^k_{\cal F}(M)=\zG_{\cal
F}(S^kTM)=\zG_{\cal F}(L_{\cal
F}M\times_{\op{GL}(n,q,\R)}S^k\R^n)$ and set \be
s(m)=[(v_1,\ldots,\linebreak v_n),(s^{i_1\ldots
i_k}(m))]_{\op{GL}(n,q,\R)},\label{s(m)}\ee where $m\in M$, and
where $(s^{i_1\ldots i_k}(m))$ is the tuple of components of
$s(m)$ in the basis induced by $(v_1,\ldots,v_n)$. Define
projection $^{\cal F}\!\zp$ by \be\begin{array}{c}(^{\cal F}\!\zp
s)(m):= [^{\cal F}\!p^1(v_1,\ldots,v_n),p_{n,q}(s^{i_1\ldots
i_k}(m))]_{\op{GL}(q,\R)}=[([v_{p+1}],\ldots,[v_{p+q}]),
(s^{p+{\frak i}_1\ldots p+{\frak i}_k}(m))]_{\op{GL}(q,\R)}\\\in
LN_m(M,{\cal
F})\times_{\op{GL}(q,\R)}S^k\R^q\end{array}\label{PiP(m)}\ee Since
$s$ is adapted, we have $^{\cal F}\!\zp s\in\zG(LN(M,{\cal
F})\times_{\op{GL}(q,\R)}S^k\R^q;{\cal F}_{S^kN})=\zG(S^kN(M,{\cal
F});{\cal F}_{S^kN})={\cal S}^k(M,{\cal F})$. Of course, $^{\cal
F}\!\zp s$ is well-defined. Indeed, if ${\cal
A}\in\op{GL}(n,q,\R)$, and if we denote the submatrices of ${\cal
A}$ by $A,B,D$, we get
$$\begin{array}{c}^{\cal F}\!p^1((v_1,\ldots,v_n)\cdot {\cal
A})=\,^{\cal F}\!p^1({\cal A}^{k_1}_1v_{k_1},\ldots,{\cal
A}^{k_n}_nv_{k_n}) =({\cal A}^{p+{\frak
k}_{p+1}}_{p+1}[v_{p+{\frak k}_{p+1}}],\ldots,{\cal A}^{p+{\frak
k}_{p+q}}_{p+q}[v_{p+{\frak k}_{p+q}}])\\=(D^{{\frak
k}_{p+1}}_{1}[v_{p+{\frak k}_{p+1}}],\ldots,D^{{\frak
k}_{p+q}}_{q}[v_{p+{\frak k}_{p+q}}])=
([v_{{p+1}}],\ldots,[v_{{p+q}}])\cdot D\end{array}$$ and \be
p_{n,q}(\zr({\cal A}^{-1})(s^{i_1\ldots i_k}(m)))=({\cal
A}^{-1,p+{\frak i}_1}_{p+{\frak j}_1} \ldots{\cal A}^{-1,p+{\frak
i}_k}_{p+{\frak j}_k}s^{p+{\frak j}_1\ldots p+{\frak
j}_k}(m))=\zr(D^{-1})(s^{p+{\frak i}_1\ldots p+{\frak
i}_k}(m)).\label{ProjActiComp}\ee

Let us recall that the isomorphism between the space
$\zG(B\times_G V)$ of sections of a vector bundle $B\times_G V\raa
M$ associated with a principal bundle $B(M,G,\zp)$ and the space
$\Ci(B,V)_G$ of $G$-equivariant functions, assigns to a section
$s$ the function $\tilde{s}$ that maps a ``basis'' $b\in B$ to the
``components'' $v$ of the ``vector'' $s(\zp(b))=[b,v]_G$ in
``basis'' $b$. It therefore follows from Equations (\ref{s(m)})
and (\ref{PiP(m)}), as well as from the definition $^{\cal
F}\!\tilde{\zp}\tilde{s}=(^{\cal F}\!\zp s)^{\tilde{}}$, that the
first part of Equation (\ref{PiFoliTildeHat}) holds true.

Eventually, $\hat{s}=\tilde{s}\circ p^2({\cal F})$ (resp.
$\hat{s}=\tilde{s}\circ p^2_{\cal F}$) in the foliated (resp.
adapted) case. The second part of Equation (\ref{PiFoliTildeHat})
is then a consequence of the definition $^{\cal
F}\!\hat{\zp}\hat{s}=(^{\cal F}\!\zp s)^{\hat{}}$, of Proposition
\ref{ProjCartBundAdapCartBundFoli}, and of the first part of
Equation (\ref{PiFoliTildeHat}).
\end{proof}

\noindent {\bf Remark}. Natural projectively invariant and
equivariant quantizations are often valued in differential
operators between tensor densities of weights $\zl$ and $\zm$.
Symbols are then sections in $\zG(STM\otimes
\zD^{\zn}TM)=\zG(LM\times_{\op{GL}(n,\R)}(S\R^n\otimes
\zD^{\zn}\R^n))$, where $\zn=\zm-\zl$, where $\zD^{\zn}TM$ is the
line bundle of $\zn$-densities on $M$ and $\zD^{\zn}\R^n$ is its
typical fiber. As the action of a change of basis, say ${\cal
A}\in\op{GL}(n,\R)$, on the component $r$ of a $\zn$-density of
$\R^n$ is, as easily checked, $\zr({\cal
A}^{-1})r=r\vert\op{det}{\cal A}\vert^{\zn}$, we get, for ${\cal
A}\in\op{GL}(n,q,\R)$,
$$p_{n,q}(\zr({\cal
A}^{-1})r)=r\vert\op{det}A\vert^{\zn}\vert\op{det}D\vert^{\zn}\neq
r\vert\op{det}D\vert^{\zn}=\zr(D^{-1})r,$$ see Equation
(\ref{ProjActiComp}). Hence, differential operators between tensor
densities, see \cite{DLO},\cite{MR}, or even between sections of
arbitrary vector bundles associated with the principal bundle of
linear frames, see \cite{BHMP},\cite{SH}, are more intricate.
Corresponding investigations are postponed to future work.

\section{Construction of the normal Cartan connection}

The method exposed in \cite{MR} in order to solve the problem of
the natural and projectively equivariant quantization uses the
notion of normal Cartan connection. We are going to adapt this
object firstly to the adapted situation and secondly to the
foliated situation. Finally, in a third step, we are going to
analyze the link between the adapted normal Cartan connection and
the foliated one.

\subsection{Construction in the adapted case}

First, recall the notion of Cartan connection on a principal fiber
bundle :

\begin{defi}
Let $G$ be a Lie group and $H$ a closed subgroup. Denote by
$\mathfrak{g}$ and $\mathfrak{h}$ the corresponding Lie algebras.
Let $P\to M$ be a principal
 $H$-bundle over $M$, such that $\mathrm{dim}\,M$ = $\mathrm{dim}\,G/H$.
 A Cartan
connection on $P$ is a $\mathfrak{g}$-valued one-form $\omega$ on
$P$ such that
\begin{itemize}
\item If $R_a$ denotes the right action of $a\in H$ on $P$, then
 $R_a^*\omega = \op{Ad}(a^{-1})\omega$,
\item If $k^*$ is the vertical vector
  field associated to $k\in\mathfrak{h}$, then $\omega (k^*)= k$,
\item $\forall u\in P, \omega_u : T_uP\to \mathfrak{g}$ is a
linear bijection.
\end{itemize}
\end{defi}

Recall too the definition of the curvature of a Cartan connection
:
\begin{defi}
If $\omega$ is a Cartan connection defined on a $H$-principal
bundle $P$, then its curvature $\Omega$ is defined as usual by
 \begin{equation}\label{curv}
\Omega = d\omega+\frac{1}{2}[\omega,\omega].
\end{equation}
\end{defi}

 \vspace{0.3cm}Next, one
adapts Theorem 4.2. cited in \cite{Koba} p.135 in the following
way :

\begin{theo}\label{Kob}
Let $P_{\mathcal{F}}$ be an $H(n+1,q+1,\mathbb{R})$-principal
fiber bundle on a manifold $M$. If one has a one-form
$\omega_{-1}$ with values in $\mathbb{R}^{n}$ of components
$\omega^{i}$ and a one-form $\omega_{0}$ with values in
$\op{gl}(n,q,\mathbb{R})$ (the Lie algebra of
$\op{GL}(n,q,\mathbb{R})$) of components $\omega_{j}^{i}$ that
satisfy the three following conditions :
\begin{itemize}
\item
$\omega_{-1}(h^{*})=0,\quad\omega_{0}(h^{*})=h_{0},\quad\forall
  h\in\op{gl}(n,q,\mathbb{R})+\mathbb{R}^{q*}$, where $h_{0}$ is the projection with respect to $\op{gl}(n,q,\mathbb{R})$ of
  $h$,
\item $(R_{a})^{*}(\omega_{-1}+\omega_{0})=(\op{Ad}\;
a^{-1})(\omega_{-1}+\omega_{0}),\quad\forall a\in
H(n+1,q+1,\mathbb{R})$, where $\op{Ad}$ $a^{-1}$ is the
application from
$\mathbb{R}^{n}+\op{gl}(n,q,\mathbb{R})+\mathbb{R}^{q*}/\mathbb{R}^{q*}$
in itself induced by the adjoint action $\op{Ad}$
  $a^{-1}$ from $\mathbb{R}^{n}+\op{gl}(n,q,\mathbb{R})+\mathbb{R}^{q*}$ into
  $\mathbb{R}^{n}+\op{gl}(n,q,\mathbb{R})+\mathbb{R}^{q*}$,
\item If $\omega_{-1}(X)=0$, then $X$ is vertical,
\end{itemize}
and the following additional condition :
\begin{equation}\label{str1}
d\omega^{i}=-\sum\omega_{k}^{i}\wedge\omega^{k},
\end{equation}
then there is a unique Cartan connection
$\omega=\omega_{-1}+\omega_{0}+\omega_{1}$ whose curvature
$\Omega$ of components $(0;\Omega_{j}^{i};\Omega_{j})$ satisfies
the following property :
$$\sum_{i=p+1}^{n}K_{jil}^{i}=0,\quad\forall j\in\{p+1,\ldots,n\},\forall l,$$
where
$$\Omega_{j}^{i}=\sum\frac{1}{2}
K_{jkl}^{i}\;\omega^{k}\wedge\omega^{l}.$$
\end{theo}
\begin{proof}
The proof goes as in \cite{Koba}. Let
$\omega=(\omega^{i};\omega_{j}^{i};\omega_{j})$ be a Cartan
connection with the given $(\omega^{i};\omega_{j}^{i})$. Thanks to
the definition of the curvature, we have
\begin{equation}\label{str2}
d\omega_{j}^{i}=-\sum\omega_{k}^{i}\wedge\omega_{j}^{k}-\omega^{i}\wedge\omega_{j}+\delta_{j}^{i}\sum\omega_{k}\wedge\omega^{k}+\Omega_{j}^{i}
\end{equation}
and
$$d\omega_{j}=-\sum\omega_{k}\wedge\omega_{j}^{k}+\Omega_{j}.$$
Applying exterior differentiation $d$ to (\ref{str1}), making use
of (\ref{str1}) and (\ref{str2}) and collecting the terms not
involving $\omega_{j}^{i}$ and $\omega_{j}$, we obtain the first
Bianchi identity :
$$\Omega_{j}^{i}\wedge\omega^{j}=0,$$
or equivalently,
$$K_{jkl}^{i}+K_{klj}^{i}+K_{lkj}^{i}=0.$$
Then the condition $\sum_{i=p+1}^{n}K_{jil}^{i}=0$ implies also
$$\sum_{i=p+1}^{n}K_{ijl}^{i}=0.$$
Now prove the uniqueness of a normal Cartan connection. Let
$\overline{\omega}=(\omega^{i};\omega_{j}^{i};\overline{\omega})$
be another Cartan connection with the given
$(\omega^{i};\omega_{j}^{i})$. Thanks to the fact that
$\overline{\omega}_{j}-\omega_{j}$ vanishes on vertical vector
fields, we can write
$$\overline{\omega}_{j}-\omega_{j}=\sum
A_{jk}\omega^{k},$$ where the coefficients $A_{jk}$ are functions
on $P$. Denoting the curvature of $\overline{\omega}$ by
$\overline{\Omega}=(0;\overline{\Omega}_{j}^{i};\overline{\Omega}_{j})$
and writing
$$\overline{\Omega}_{j}^{i}=\sum\frac{1}{2}
\overline{K}_{jkl}^{i}\;\omega^{k}\wedge\omega^{l},$$ we obtain
using (\ref{str2}) the following relations between $K_{jkl}^{i}$
and $\overline{K}_{jkl}^{i}$ :
$$\overline{K}_{jkl}^{i}-K_{jkl}^{i}=-\delta_{l}^{i}A_{jk}+\delta_{k}^{i}A_{jl}+\delta_{j}^{i}A_{kl}-\delta_{j}^{i}A_{lk}.$$
Hence,
\begin{equation}\label{equ1}
\sum_{i=p+1}^{n}(\overline{K}_{ikl}^{i}-K_{ikl}^{i})=(q+1)(A_{kl}-A_{lk}),
\end{equation}
\begin{equation}\label{equ2}
\sum_{i=p+1}^{n}(\overline{K}_{jil}^{i}-K_{jil}^{i})=(q-1)A_{jl}+(A_{jl}-A_{lj}).
\end{equation}
If $\omega$ and $\overline{\omega}$ are normal Cartan connections,
i.e.,
$\sum_{i=p+1}^{n}\overline{K}_{jil}^{i}=\sum_{i=p+1}^{n}K_{jil}^{i}=0$,
then $A_{ij}=0$ and hence $\omega=\overline{\omega}$. This prove
the uniqueness of the normal Cartan connection.

To prove the existence, one assumes that there is a Cartan
connection $\omega=(\omega^{i};\omega_{j}^{i};\omega_{j})$ with
the given $(\omega^{i};\omega_{j}^{i})$. The goal is then to find
functions $A_{jk}$ such that
$\overline{\omega}=(\omega^{i};\omega_{j}^{i};\overline{\omega}_{j})$
becomes a normal Cartan connection. If $1\leq j\leq p$, $A_{jk}$
is of course equal to zero. If $p+1\leq j\leq n$ and if $p+1\leq
k\leq n$, one can view thanks to (\ref{equ1}) and (\ref{equ2})
that it suffices to set
\begin{equation}\label{equ3}
A_{jk}=\frac {1}{(q+1)(q-1)}\sum_{i=p+1}^{n}K_{ijk}^{i}-\frac
{1}{q-1}\sum_{i=p+1}^{n}K_{jik}^{i}.
\end{equation}
If $1\leq k\leq p$, one sees thanks to (\ref{equ2}) that it
suffices to set
\begin{equation}\label{equ4}
A_{jk}=-\sum_{i=p+1}^{n}\frac {1}{q}K_{jik}^{i}.
\end{equation}
The last step of the proof consists in showing that there is at
least one Cartan connection $\omega$ with the given
$(\omega^{i},\omega_{j}^{i})$. Let $\{U_{\alpha}\}$ be a locally
finite open cover of $M$ with a partition of unity
$\{f_{\alpha}\}$. If $\omega_{\alpha}$ is a Cartan connection in
$P_{\mathcal{F}}|U_{\alpha}$ with the given
$(\omega^{i};\omega_{j}^{i})$, then
$\sum_{\alpha}(f_{\alpha}\circ\pi)\omega_{\alpha}$ is a Cartan
connection in $P_{\mathcal{F}}$ with the given
$(\omega^{i};\omega_{j}^{i})$, where
$\pi:P_{\mathcal{F}}\rightarrow M$ is the projection. Hence, the
problem is reduced to the case where $P_{\mathcal{F}}$ is a
product bundle. Fixing a cross section $\sigma:M\rightarrow
P_{\mathcal{F}}$, set $\omega_{j}(X)=0$ for every vector $X$
tangent to $\sigma(M)$. If $Y$ is an arbitrary tangent vector of
$P_{\mathcal{F}}$, we can write uniquely
$$Y=R_{a}(X)+W,$$
where $X$ is a vector tangent to $\sigma(M)$, $a$ is in
$H(n+1,q+1,\mathbb{R})$ and $W$ is a vertical vector. Extend $W$
to a unique fundamental vector field $A^{*}$ of $P_{\mathcal{F}}$
with $A\in\op{gl}(n,q,\mathbb{R})+\mathbb{R}^{q*}$. Thanks to the
properties of the Cartan connections, we have to set
$$\omega(Y)=\op{Ad}(a^{-1})(\omega(X))+A.$$
This defines the desired $(\omega_{j})$. Actually, $X$ is equal to
$(\sigma_{*}\circ\pi_{*})Y$ and one can take the section $\sigma$
equal to $\sigma_{\alpha}$.
\end{proof}
One can remark that the codimension of the foliation $\mathcal{F}$
has to be different from 1.

One can define on $P_{\mathcal{F}}$ an one-form in the following
way :
\begin{defi}
If $u=j_{0}^{2}f$ is a point belonging to $P_{\mathcal{F}}$ and if
$X$ is a tangent vector to $P_{\mathcal{F}}$ at $u$, the canonical
form $\theta_{\mathcal{F}}$ of $P_{\mathcal{F}}$ is the 1-form
with values in $\mathbb{R}^{n}\oplus\op{gl}(n,q,\mathbb{R})$
defined at the point $u$ in the following way :
$$\theta_{\mathcal{F};u}(X)=(P^{1}f)_{*e}^{-1}(p_{\mathcal{F}*}^{2}X),$$
where $e$ is the frame at the origin of $\R^{n}$ represented by
the identity matrix.
\end{defi}
\begin{theo}
One can associate to the projective class of an adapted connection
$[\nabla_{\mathcal{F}}]$ a Cartan connection on $P_{\mathcal{F}}$
in a natural way. We will denote by $\omega_{\mathcal{F}}$ this
Cartan connection.
\end{theo}
\begin{proof}
The canonical one-form defined above is the restriction to
$P_{\mathcal{F}}$ of the canonical one-form of $P^{2}(M)$ defined
in \cite{Koba} p.140. It is too the restriction to
$P_{\mathcal{F}}$ of the restriction to $P$ of the canonical
one-form of $P^{2}(M)$, where $P$ is the projective structure
associated to $\nabla_{\mathcal{F}}$ defined in \cite{Koba}.
Thanks to the fact that the canonical one-form on $P$ satisfies
the properties of Theorem 4.2. mentioned in \cite{Koba},
$\theta_{\mathcal{F}}$ satisfies the properties mentioned in
Theorem \ref{Kob}. One defines then the adapted normal Cartan
connection $\omega_{\mathcal{F}}$ as the unique Cartan connection
on $P_{\mathcal{F}}$ beginning by $\theta_{\mathcal{F}}$ and
satisfying the property linked to the curvature cited in Theorem
\ref{Kob}. Because of the naturality of this property, the
naturality of $\theta_{\mathcal{F}}$ and the uniqueness of the
Cartan connection mentioned in Theorem \ref{Kob},
$\omega_{\mathcal{F}}$ is a Cartan connection on $P_{\mathcal{F}}$
associated naturally to the class $[\nabla_{\mathcal{F}}]$.
\end{proof}
\subsection{Construction in the foliated case}

 \vspace{0.2cm}The reduction $P(\mathcal{F})$ is actually an
 example of a foliated bundle defined in \cite{Blum}. The Cartan
 connection that we are going to define on it is an example of a
 Cartan connection in a foliated bundle defined too in
 \cite{Blum}. It is the reason for which we are going first to recall
 the definitions of these notions.

\begin{defi}Let $M$ be a manifold of dimension $m$ and let $\mathcal{F}$
be a codimension $q$ foliation of $M$. Let $T(M)$ be the tangent
bundle of $M$ and let $T\mathcal{F}$ be the tangent bundle of
$\mathcal{F}$. Let $H$ be a Lie group and let $\pi:P\to M$ be a
principal $H$-bundle. We say $\pi:P\to M$ is a foliated bundle if
there is a foliation $\tilde{\mathcal{F}}$ of $P$ satisfying
\begin{itemize}
\item $\tilde{\mathcal{F}}$ is $H$-invariant, \item
$\tilde{E}_{u}\cap V_{u}=\{0\}$ for all $u\in P$, \item
$\pi_{*u}(\tilde{E}_{u})=T\mathcal{F}_{\pi(u)}$ for all $u\in P$,
\end{itemize}
where $\tilde{E}$ is the tangent bundle of $\tilde{\mathcal{F}}$
and $V$ is the bundle of vertical vectors.
\end{defi}
\begin{defi}
Let $\mathcal{F}$ be a codimension $q$ foliation of $M$. Let $G$
be a Lie group and let $H$ be a closed subgroup of $G$ with
dimension($G/H)=q$. Let $\pi:P\to M$ be a foliated principal
$H$-bundle. Let $\mathfrak{g}$ be the Lie algebra of $G$ and let
$\mathfrak{h}$ be the Lie algebra of $H$. For each
$A\in\mathfrak{h}$, let $A^{*}$ be the corresponding fundamental
vector field on $P$.

A Cartan connection in the foliated bundle $\pi:P\to M$ is a
$\mathfrak{g}$-valued one-form $\omega$ on $P$ satisfying
\begin{itemize}
\item $\omega(A^{*})=A$ for all $A\in\mathfrak{h}$, \item
$(R_{a})^{*}\omega=\op{Ad}(a^{-1})\omega$ for all $a\in H$ where
$R_{a}$ denotes the right translation by $a$ acting on $P$ and
$\op{Ad}(a^{-1})$ is the adjoint action of $a^{-1}$ on
$\mathfrak{g}$, \item For each $u\in P$,
$\omega_{u}:T_{u}P\to\mathfrak{g}$ is onto and
$\omega_{u}(\tilde{E}_{u})=0$, \item $L_{X}\omega=0$ for all
$X\in\Gamma(\tilde{E})$ where $\Gamma(\tilde{E})$ denotes the
smooth sections of $\tilde{E}$ and $L_{X}$ is the Lie derivative.
\end{itemize}
\end{defi}
\begin{theo}
The reduction $P(\mathcal{F})$ is a foliated bundle.
\end{theo}
\begin{proof}
One can easily view that $\mathcal{F}_{P^{2}N}$ satisfies the
properties of the definition of a foliated bundle : first,
$\mathcal{F}_{P^{2}N}$ is $H(q+1,\mathbb{R})$-invariant because,
if $(U_{i},f_{i},g_{ij})$ is a cocycle corresponding to the
foliation $\mathcal{F}$, if $J_{0}^{2}f\in P(\mathcal{F})$ and if
$j_{0}^{2}(f_i\circ f)$ is constant, then $j_{0}^{2}(f_i\circ
f\circ h)=j_{0}^{2}(f_i\circ f)\circ j_{0}^{2}(h)$ is constant,
where $j_{0}^{2}(h)\in H(q+1,\mathbb{R})$.

\vspace{0.2cm}If $X\in V_{u}$, then $X=\frac
{d}{dt}u\exp(th)|_{t=0}$, where $h\in
\mathfrak{h}(q+1,\mathbb{R})$. If $u=J_{0}^{2}(f)$ and if
$\exp(th)=j_{0}^{2}(g_t)$, then $j_{0}^{2}(f_i\circ f\circ g_{t})$
is constant if $X$ is tangent to the foliation
$\mathcal{F}_{P^{2}N}$. One has then that $J_{0}^{2}(f\circ
g_{t})$ is constant and then $X=0$.

\vspace{0.2cm}If $X$ is tangent to the foliation
$\mathcal{F}_{P^{2}N}$, then $X=\frac {d}{dt}\gamma(t)|_{t=0}$,
where $\gamma(t)\in\mathcal{F}_{P^{2}N}$. Then
$\pi_{*u}^{2}(X)=\frac {d}{dt}\pi^{2}(\gamma(t))|_{t=0}$, that
belongs to $T\mathcal{F}_{\pi(u)}$ because $\pi^{2}(\gamma(t))$
belongs to $\mathcal{F}$. Indeed, if $\gamma(t)=J_{0}^{2}(f_t)$,
$f_i\circ\pi^{2}(\gamma(t))$ is constant because
$j_{0}^{2}(f_i\circ f_{t})$ is constant.
\end{proof}
\begin{theo}
One can associate to the class of a foliated connection
$[\nabla(\mathcal{F})]$ a Cartan connection on $P(\mathcal{F})$ in
a natural way. We will denote this connection by
$\omega(\mathcal{F})$.
\end{theo}
\begin{proof}
If $(U_{i},f_{i},g_{ij})$ is a Haefliger cocycle of $\mathcal{F}$,
the image by $P^{2}N(f_i)$ of $P(\mathcal{F})|_{U_{i}}$ is a
reduction of $P^{2}N(f_i)(P^{2}N(U_i,\mathcal{F}))$ to
$H(q+1,\mathbb{R})$. Indeed, if $j_{0}^{2}(h)\in
H(q+1,\mathbb{R})$, $j_{0}^{2}(f_i\circ
f)j_{0}^{2}(h)=j_{0}^{2}(f_i\circ f\circ
h)=P^{2}N(f_i)(J_{0}^{2}f\circ j_{0}^{2}(h))$, with
$J_{0}^{2}f\circ j_{0}^{2}(h)\in P(\mathcal{F})$ if $J_{0}^{2}f\in
P(\mathcal{F})$. Moreover, if $j_{0}^{2}(f_i\circ f')$ and
$j_{0}^{2}(f_i\circ f)$ belong to the same fiber, then $f_i\circ
f'(0)=f_i\circ f(0)$. If $y$ denotes the passing to the transverse
coordinates of an adapted coordinates system, one has then
$j_{0}^{2}(y\circ f')=(x^{i},\delta_{k}^{i},-\Gamma_{jk}^{i})H'$
and $ j_{0}^{2}(y\circ
f)=(x^{i},\delta_{k}^{i},-\Gamma_{jk}^{i})H$, with $H$ and $H'$
belonging to $H(q+1,\mathbb{R})$ and with the $\Gamma_{jk}^{i}$
equal to the Christoffel symbols of $\nabla(\mathcal{F})$. We have
thus $ j_{0}^{2}(f_i\circ f')= j_{0}^{2}(f_i\circ f)H^{-1}H'$,
with $H^{-1}H'\in H(q+1,\mathbb{R})$.

\vspace{0.2cm}We will denote by $\overline{P}$ the reduction of
$P^{2}N(f_i)(P^{2}N(U_i,\mathcal{F}))$ to $H(q+1,\mathbb{R})$.

\vspace{0.2cm}One builds locally the normal Cartan connection
$\omega(\mathcal{F})$ on $P(\mathcal{F})$ in the following way :
if $\overline{\omega}$ denotes the normal Cartan connection on
$\overline{P}$, then
$\omega(\mathcal{F})|_{P^{2}N(U_i,\mathcal{F})}:=(P^{2}N(f_i))^{*}\overline{\omega}$.

\vspace{0.2cm}One can show (see \cite{Blum}) that the connection
$\omega(\mathcal{F})$ is a well-defined foliated Cartan
connection.

\vspace{0.2cm}Thanks to the naturality of the normal Cartan
connection, $\omega(\mathcal{F})$ is associated naturally to the
class of the foliated connection $[\nabla(\mathcal{F})]$.
\end{proof}
 \vspace{0.2cm}One can remark that, as the
foliation $\mathcal{F}_{P^{2}N}$ is of dimension $p$, the third
condition of the definition of a foliated Cartan connection
implies that, in our case, the kernel of $\omega(\mathcal{F})_{u}$
will be exactly equal to the tangent space to
$\mathcal{F}_{P^{2}N}$.
\subsection{Link between adapted and foliated Cartan connections}
\noindent {\bf Remark}. If $y$ denotes the passing to the
transverse coordinates of an adapted coordinates system and if
$P^{2}y$ denotes the following application :
$$P^{2}y:P^{2}(M,\mathcal{F})\to P^{2}\mathbb{R}^{q}:J_{0}^{2}f\mapsto j_{0}^{2}(y\circ f),$$
 the image by $P^{2}y$ of
$P(\mathcal{F})$ is a reduction of $P^{2}(U)$ to
$H(q+1,\mathbb{R})$, where $U$ is an open set of $\mathbb{R}^{q}$.
We will denote by $P_{U}$ this reduction of $P^{2}(U)$ to
$H(q+1,\mathbb{R})$. If $\omega_{U}$ denotes the normal Cartan
connection on $P_{U}$, then
$\omega(\mathcal{F})_{(P^{2}y)^{-1}P_{U}}=(P^{2}y)^{*}\omega_{U}$.
Indeed, if $\phi$ denotes the diffeomorphism such that $\phi\circ
y=f_{i}$, then $P^{2}\phi(P_{U})=\overline{P}$. By naturality of
the normal Cartan connection,
$\omega_{U}=(P^{2}\phi)^{*}\overline{\omega}$ and then
$(P^{2}y)^{*}\omega_{U}=(P^{2}f_i)^{*}\overline{\omega}$.
\begin{prop}\label{lemme}
If $\theta_{U}$ denotes the canonical one-form on $P_{U}$, then
$$(P^{2}y\circ(^{\cal
F}\!\!{\frak p}^2))^{*}\theta_{U}=p_{n,q}\theta_{\mathcal{F}}.$$
\end{prop}
\begin{proof}
On one hand, if $u=j_{0}^{2}(f_0)\in P_{\mathcal{F}}$ and if
$X=\frac {d}{dt}j_{0}^{2}(f_t)|_{t=0}\in T_{u}P_{\mathcal{F}}$,
one has $\theta_{\mathcal{F}\;u}(X)=
((P^{1}f_0)^{-1})_{*}(p_{\mathcal{F}*}^{2}X)=[(P^{1}f_0)^{-1}\circ
p_{\mathcal{F}}^{2}]_{*u}\frac {d}{dt}j_{0}^{2}(f_t)|_{t=0}=\frac
{d}{dt}(P^{1}f_0)^{-1}j_{0}^{1}(f_t)|_{t=0}=\frac
{d}{dt}j_{0}^{1}(f_0^{-1}\circ f_{t})|_{t=0}$.

\vspace{0.2cm}On the other hand, $((P^{2}y\circ(^{\cal
F}\!\!{\frak
p}^2))^{*}\theta_{U})_{u}(X)=\theta_{U\,P^{2}y\circ(^{\cal
F}\!\!{\frak p}^2)(u)}(P^{2}y_{*}{^{\cal F}\!\!{\frak p}^2}_{*}X)$
is equal to $\frac {d}{dt}(P^{1}(y\circ f_0\circ i_{q})^{-1}\circ
p_{U}^{2}\circ P^{2}y\circ(^{\cal F}\!\!{\frak
p}^2)(j_{0}^{2}(f_t)))|_{t=0}=\frac {d}{dt}j_{0}^{1}((y\circ
f_0\circ i_{q})^{-1}\circ(y\circ f_{t}\circ i_{q}))|_{t=0}$, if
$p_{U}^{2}$ denotes the projection of $P_{U}$ on $P^{1}(U)$. One
can then easily show that $p_{n,q}\frac
{d}{dt}f_0^{-1}(f_t(0))|_{t=0}=\frac {d}{dt}(y\circ f_0\circ
i_{q})^{-1}y(f_t(0))|_{t=0}$ and that $p_{n,q}\frac
{d}{dt}[(f_0^{-1}\circ f_{t})_{*0}]|_{t=0}=\frac {d}{dt}[((y\circ
f_0\circ i_q)^{-1}\circ(y\circ f_{t}\circ i_q))_{*0}]|_{t=0}$
using the fact that the differentials of the local forms of $f_0$
and $f_{t}$ belong to $\op{GL}(n,q,\mathbb{R})$.
\end{proof}
\begin{theo}
The connections $\omega_{\mathcal{F}}$ and $\omega(\mathcal{F})$
are linked by the following relation :
$${^{\cal F}\!\!{\frak p}^2}^{*}\omega(\mathcal{F})=p_{n,q}\omega_{\mathcal{F}}.$$
\end{theo}
\begin{proof}
To prove that, it suffices to prove that
$$(P^{2}y\circ(^{\cal
F}\!\!{\frak p}^2))^{*}\omega_{U}=p_{n,q}\omega_{\mathcal{F}}.$$

If one denotes by $\nabla_{U}$ the connection on $U$ whose
Christoffel symbols are the Christoffel symbols of
$\nabla(\mathcal{F})$ (the $\Gamma_{jk}^{i}$ with $i,j,k$ between
$p+1$ and $n$), if $(\epsilon^{1},\ldots,\epsilon^{n})$ denotes
the canonical basis of $\mathbb{R}^{n*}$ (resp.
$(\epsilon^{1},\ldots,\epsilon^{q})$ denotes the canonical basis
of $\mathbb{R}^{q*}$), one has
$$\omega_{\mathcal{F}}=\tilde{\Upsilon}_{\mathcal{F}}-\sum_{j=p+1}^{n}\sum_{k=1}^{n}(\Gamma_{\mathcal{F}\;jk})(\theta_{\mathcal{F}-1}^{k})\epsilon^{j}$$
$$(resp.\;\omega_{U}=\tilde{\Upsilon}_{U}-\sum_{j=1}^{q}\sum_{k=1}^{q}(\Gamma_{U\;jk})(\theta_{U\;-1}^{k})\epsilon^{j}),$$

where $\tilde{\Upsilon}_{\mathcal{F}}$ (resp.
$\tilde{\Upsilon}_{U}$) is the Cartan connection induced by
$\nabla_{\mathcal{F}}$ (resp. $\nabla_{U}$),
$\Gamma_{\mathcal{F}}$ (resp. $\Gamma_{U}$) is the deformation
tensor corresponding to $\nabla_{\mathcal{F}}$ (resp.
$\nabla_{U}$) (see \cite{Capinv}).

\vspace{0.2cm}One recall that $\tilde{\Upsilon}_{\mathcal{F}}$
(resp. $\tilde{\Upsilon}_{U}$) is the unique Cartan connection
such that its component with respect to $\mathbb{R}^{q*}$ vanishes
on the section $(x^{i},\delta_{k}^{i},-\Gamma_{jk}^{i}$). If
$\sigma_{\mathcal{F}}$ (resp. $\sigma_{U}$) denotes the section
$(x^{i},\delta_{k}^{i},-\Gamma_{jk}^{i})$, the connection
$\tilde{\Upsilon}_{\mathcal{F}}$ (resp. $\tilde{\Upsilon}_{U}$) is
defined in this way :
$$\tilde{\Upsilon}_{\mathcal{F}\;u}(X)=\op{Ad}(b^{-1})\theta_{\mathcal{F}}((\sigma_{\mathcal{F}}\circ\pi^{2})_{*}X)+B,$$
$$(resp.\;\tilde{\Upsilon}_{U\;u}(X)=\op{Ad}(b^{-1})\theta_{U}((\sigma_{U}\circ\pi^{2})_{*}X)+B),$$

where $\pi^{2}$ is the projection on $M$ (resp. $U$),
$R_{b}(\sigma_{\mathcal{F}}(\pi^{2}(u)))=u$ (resp.
$R_{b}(\sigma_{U}(\pi^{2}(u)))=u$) and

\vspace{0.2cm}$B^{*}=X-R_{b*}\sigma_{\mathcal{F}*}\pi_{*}^{2}X$
(resp. $B^{*}=X-R_{b*}\sigma_{U*}\pi_{*}^{2}X$).

\vspace{0.2cm}The deformation tensor $\Gamma_{\mathcal{F}}$ (resp.
$\Gamma_{U}$) is defined in this way :
$$\Gamma_{\mathcal{F}}(X)=(\tilde{\Upsilon}_{\mathcal{F}}-\omega_{\mathcal{F}})(\omega_{\mathcal{F}}^{-1}(X))$$
$$(resp.\;\Gamma_{U}(X)=(\tilde{\Upsilon}_{U}-\omega_{U})(\omega_{U}^{-1}(X))).$$

\vspace{0.2cm}In fact, the sections
$(x^{i},\delta_{k}^{i},-\Gamma_{jk}^{i})$ correspond to the
section $\sigma$ of the end of the theorem \ref{Kob}, the
connections $\tilde{\Upsilon}_{\mathcal{F}}$ and
$\tilde{\Upsilon}_{U}$ correspond to the connection $\omega$ of
the proof of this theorem, the connections $\omega_{\mathcal{F}}$
and $\omega_{U}$ correspond to the connection $\bar{\omega}$
whereas the $\Gamma_{jk}$ correspond to the functions $-A_{jk}$.

 \vspace{0.2cm}We first prove that $(P^{2}y\circ(^{\cal F}\!\!{\frak p}^2))^{*}\tilde{\Upsilon}_{U}=p_{n,q}\tilde{\Upsilon}_{\mathcal{F}}$.

\vspace{0.2cm}Indeed, we have
$$\tilde{\Upsilon}_{\mathcal{F}\;u}(X)=\op{Ad}(b^{-1})\theta_{\mathcal{F}}((\sigma_{\mathcal{F}}\circ\pi^{2})_{*}X)+B,$$
where $R_{b}(\sigma_{\mathcal{F}}(\pi^{2}(u)))=u$ and
$B^{*}=X-R_{b*}\sigma_{\mathcal{F}*}\pi_{*}^{2}X$ whereas
$$(P^{2}y\circ(^{\cal F}\!\!{\frak p}^2))^{*}\tilde{\Upsilon}_{U\;u}(X)=\op{Ad}(a^{-1})\theta_{U}((\sigma_{U}\circ\pi^{2})_{*}(P^{2}y\circ(^{\cal F}\!\!{\frak p}^2))_{*}X)+A,$$
where $R_{a}(\sigma_{U}(\pi^{2}(P^{2}y\circ(^{\cal F}\!\!{\frak
p}^2))(u)))=(P^{2}y\circ(^{\cal F}\!\!{\frak p}^2))u$ and
$A^{*}=(P^{2}y\circ(^{\cal F}\!\!{\frak
p}^2))_{*}X-R_{a*}\sigma_{U*}\pi_{*}^{2}(P^{2}y\circ(^{\cal
F}\!\!{\frak p}^2))_{*}X$.

\vspace{0.2cm}One can see that $a=p_{n,q}b$.

\vspace{0.2cm}Moreover, as
$(\sigma_{U}\circ\pi^{2})\circ(P^{2}y\circ(^{\cal F}\!\!{\frak
p}^2))=(P^{2}y\circ(^{\cal F}\!\!{\frak
p}^2))\circ(\sigma_{\mathcal{F}}\circ\pi^{2})$ and
$(P^{2}y\circ(^{\cal F}\!\!{\frak
p}^2))^{*}\theta_{U}=p_{n,q}\theta_{\mathcal{F}}$ (see proposition
\ref{lemme}), one has
$\op{Ad}((p_{n,q}b)^{-1})\theta_{U}((\sigma_{U}\circ\pi^{2})_{*}(P^{2}y\circ(^{\cal
F}\!\!{\frak
p}^2))_{*}X)=\op{Ad}((p_{n,q}b)^{-1})(p_{n,q}\theta_{\mathcal{F}})((\sigma_{\mathcal{F}}\circ\pi^{2})_{*}X)=p_{n,q}(\op{Ad}(b^{-1})\theta_{\mathcal{F}}((\sigma_{\mathcal{F}}\circ\pi^{2})_{*}X))$.
One can see too that $A^{*}=(P^{2}y\circ(^{\cal F}\!\!{\frak
p}^2))_{*}B^{*}$, thus $A=p_{n,q}B$.

\vspace{0.2cm}Now, prove that $(P^{2}y\circ(^{\cal F}\!\!{\frak
p}^2))^{*}\sum_{j=1}^{q}\sum_{k=1}^{q}(\Gamma_{U\;jk})(\theta_{U\;-1}^{k})\epsilon^{j}=\sum_{j=p+1}^{n}\sum_{k=1}^{n}(\Gamma_{\mathcal{F}\;jk})(\theta_{\mathcal{F}\;-1}^{k})\epsilon^{j}$.

\vspace{0.2cm}One has $(P^{2}y\circ(^{\cal F}\!\!{\frak
p}^2))^{*}\sum_{j=1}^{q}\sum_{k=1}^{q}(\Gamma_{U\;jk})(\theta_{U\;-1}^{k})\epsilon^{j}=\sum_{j=1}^{q}\sum_{k=1}^{q}(P^{2}y\circ(^{\cal
F}\!\!{\frak
p}^2))^{*}(\Gamma_{U\;jk})(\theta_{\mathcal{F}\;-1}^{k+p})\epsilon^{j+p}$.

\vspace{0.2cm}It remains then to prove that $(P^{2}y\circ(^{\cal
F}\!\!{\frak
p}^2))^{*}(\Gamma_{U\;jk})=\Gamma_{\mathcal{F}\;j+p,k+p}$ and that
$\Gamma_{\mathcal{F}\;jk}=0$ if $1\leq k\leq p$.

\vspace{0.2cm}Indeed, if $1\leq k\leq p$,
$\Gamma_{\mathcal{F}\;jk}=\frac
{1}{q}\sum_{l=p+1}^{n}R_{\mathcal{F}\,jlk}^{l}$, where
$R_{\mathcal{F}}$ denotes the equivariant function on
$P_{\mathcal{F}}$ representing the curvature tensor of
$\nabla_{\mathcal{F}}$ thanks to the equation (\ref{equ4}) of the
Theorem \ref{Kob} and thanks to the fact that the $K_{ijk}^{l}$
represent the components of $R_{\mathcal{F}}$ (see \cite{Capinv}).
Thanks to the fact that $\nabla_{\mathcal{F}}$ is adapted, one can
see that if $1\leq k\leq p$, $\Gamma_{\mathcal{F}\;jk}=0$.
Moreover, $\Gamma_{U\;jk}=\frac
{-1}{(q+1)(q-1)}\sum_{i=1}^{q}R_{U\,ijk}^{i}+\frac
{1}{q-1}\sum_{i=1}^{q}R_{U\,jik}^{i}$, where $R_{U}$ denotes the
equivariant function on $P_{U}$ representing the curvature tensor
of $\nabla_{U}$ whereas if $p+1\leq k\leq n$,
$\Gamma_{\mathcal{F}\;jk}=\frac
{-1}{(q+1)(q-1)}\sum_{i=p+1}^{n}R_{\mathcal{F}\,ijk}^{i}+\frac
{1}{q-1}\sum_{i=p+1}^{n}R_{\mathcal{F}\,jik}^{i}$ thanks to the
equation (\ref{equ3}) of the Theorem \ref{Kob}. This allows to
prove that $(P^{2}y\circ(^{\cal F}\!\!{\frak
p}^2))^{*}(\Gamma_{U\;jk})=\Gamma_{\mathcal{F}\;j+p,k+p}$.

\end{proof}
 \section{Construction of
the quantization} In a first step, we are going to explain how to
build the quantization in the adapted and foliated situations. In
a second step, we are going to prove that the quantization
commutes with the reduction. In other words, quantize adapted
objects is equivalent to quantize the induced foliated objects.
\subsection{Construction in the adapted situation}
In the adapted situation, we can define the operator of invariant
differentiation exactly in the same way as in the standard
situation :
\begin{defi}
Let $V$ be a vector space. If $f\in
 C^{\infty}(P_{\mathcal{F}},V)$, then the invariant differential of $f$ with respect to
$\omega_{\mathcal{F}}$ is the function
$\nabla^{\omega_{\mathcal{F}}}f\in
C^{\infty}(P_{\mathcal{F}},\mathbb{R}^{n*}\otimes V)$ defined by
\[ \nabla^{\omega_{\mathcal{F}}}f(u)(X) = L_{\omega_{\mathcal{F}}^{-1}(X)}f(u)\quad\forall u\in
P_{\mathcal{F}},\quad\forall X\in\R^n.\]
\end{defi}
We will also use an iterated and symmetrized version of the
invariant differentiation
\begin{defi}
If $f\in C^{\infty}(P_{\mathcal{F}},V)$ then
$(\nabla^{\omega_{\mathcal{F}}})^k f \in
C^{\infty}(P_{\mathcal{F}},S^k\R^{n*}\otimes V)$ is defined by
\[(\nabla^{\omega_{\mathcal{F}}})^k f(u)(X_1,\ldots,X_k) = \frac{1}{k!}\sum_{\nu}
L_{\omega_{\mathcal{F}}^{-1}(X_{\nu_1})}\circ\ldots\circ
L_{\omega_{\mathcal{F}}^{-1}(X_{\nu_k})}f(u)\] for
$X_1,\ldots,X_k\in\R^n$.
\end{defi}

\begin{prop}
If $v\in\mathbb{R}^{n}$ and if $p_{n,q}(v)=0$, then
$\omega_{\mathcal{F}}^{-1}(v)$ is tangent to
$\mathcal{F}_{\underline{P}^{2}}$. \end{prop}
\begin{proof}
 Indeed, as
$p_{n,q}\omega_{\mathcal{F}}={^{\cal F}\!\!{\frak
p}^2}^{*}\omega(\mathcal{F})$, one has
$\omega(\mathcal{F})({^{\cal F}\!\!{\frak
p}^2}_{*}\omega_{\mathcal{F}}^{-1}(v))=0$. As ${^{\cal
F}\!\!{\frak p}^2}_{*}\omega_{\mathcal{F}}^{-1}(v)$ is then
tangent to $\mathcal{F}_{P^{2}N}$, one can easily show that
$\omega_{\mathcal{F}}^{-1}(v)$ is then tangent to
$\mathcal{F}_{\underline{P}^{2}}$.
\end{proof}
In the adapted situation, the invariant differentiation has a
particular property :
\begin{prop}\label{inv}
If $f$ is a foliated function on $P_{\mathcal{F}}$, then
$$(\nabla^{\omega_{\mathcal{F}}^{k}}f)(v_{1},\ldots,v_{k})=(\nabla^{\omega_{\mathcal{F}}^{k}}f)((0,p_{n,q}v_{1}),\ldots,(0,p_{n,q}v_{k})).$$
\end{prop}
\begin{proof}
Indeed, one can show that if $f$ is constant along the leaves of
$\mathcal{F}_{\underline{P}^{2}}$, then
$L_{\omega_{\mathcal{F}}^{-1}(0,p_{n,q}v)}f$ is a foliated
function too if $v\in\mathbb{R}^{n}$. Indeed, if $X$ is tangent to
$\mathcal{F}_{\underline{P}^{2}}$, then
$L_{X}L_{\omega_{\mathcal{F}}^{-1}(0,p_{n,q}v)}f=0$. To show that,
it suffices to prove that
$L_{[X,\omega_{\mathcal{F}}^{-1}(0,p_{n,q}v)]}f=0$. The fact that
$i_{X}\omega(\mathcal{F})=i_{X}d\omega(\mathcal{F})=0$ if $X$ is
tangent to $\mathcal{F}_{P^{2}N}$, that
$p_{n,q}\omega_{\mathcal{F}}={^{\cal F}\!\!{\frak
p}^2}^{*}\omega(\mathcal{F})$ and that ${^{\cal F}\!\!{\frak
p}^2}_{*}X$ is tangent to $\mathcal{F}_{P^{2}N}$ if $X$ is tangent
to $\mathcal{F}_{\underline{P}^{2}}$ implies that
$i_{X}p_{n,q}\omega_{\mathcal{F}}=i_{X}p_{n,q}d\omega_{\mathcal{F}}=0$
if $X$ is tangent to $\mathcal{F}_{\underline{P}^{2}}$.

\vspace{0.2cm}Remark that as the kernel of
$p_{n,q}\omega_{\mathcal{F}}$ has a dimension equal to the
dimension of $\mathcal{F}_{\underline{P}^{2}}$ (i.e. $p+np$), the
kernel of $p_{n,q}\omega_{\mathcal{F}}$ is equal to the tangent
space to $\mathcal{F}_{\underline{P}^{2}}$. One has then
$0=p_{n,q}d\omega_{\mathcal{F}}(X,\omega_{\mathcal{F}}^{-1}(0,p_{n,q}v))=X.(p_{n,q}v)-\omega_{\mathcal{F}}^{-1}(0,p_{n,q}v).(p_{n,q}\omega_{\mathcal{F}}(X))-p_{n,q}\omega_{\mathcal{F}}([X,\omega_{\mathcal{F}}^{-1}(0,p_{n,q}v)])$.
As the first two terms are equal to 0, the third term vanishes
too.

\vspace{0.2cm}One has then that
$[X,\omega_{\mathcal{F}}^{-1}(0,p_{n,q}v)]$ is tangent to
$\mathcal{F}_{\underline{P}^{2}}$ and then
$L_{[X,\omega_{\mathcal{F}}^{-1}(0,p_{n,q}v)]}f=0$.

\vspace{0.2cm}One concludes using the fact that
$\omega_{\mathcal{F}}^{-1}(v)$ is tangent to
$\mathcal{F}_{\underline{P}^{2}}$ if $p_{n,q}(v)=0$.
\end{proof}

 \vspace{0.2cm}In the adapted situation, we define a divergence
operator analogous to the divergence operator defined in
\cite{MR}.

We fix a basis $(e_1,\ldots,e_n)$ of $\R^n$ and we denote by
$(\epsilon^1,\ldots,\epsilon^n)$ the dual basis in $\R^{n*}$.
\begin{defi}
The \emph{Divergence operator} with respect to the Cartan
connection $\omega_{\mathcal{F}}$ is defined by
\[\op{Div}^{\omega_{\mathcal{F}}} : C^{\infty}(P_{\mathcal{F}},S^k(\R^n))
\to C^{\infty}(P_{\mathcal{F}},S^{k-1}(\R^n)) : S\mapsto
\sum_{j=p+1}^n
i(\epsilon^j)\nabla^{\omega_{\mathcal{F}}}_{e_j}S,\] where $i$
denotes the inner product.
\end{defi}

\noindent {\bf Remark}. If $S\in
C^{\infty}(P_{\mathcal{F}},S^k(\mathbb{R}^n))$ and if $f\in
C^{\infty}(P_{\mathcal{F}},\mathbb{R};\mathcal{F}_{\underline{P}^{2}})$,
thanks to Proposition \ref{inv}, we have
$\langle\op{Div}^{\omega_{\mathcal{F}}^{l}}S,
\nabla^{\omega_{\mathcal{F}}^{k-l}}f\rangle=\langle
p_{n,q}\op{Div}^{\omega_{\mathcal{F}}^{l}}S,
p_{n,q}\nabla^{\omega_{\mathcal{F}}^{k-l}}f\rangle$.

\vspace{0.2cm}One can then easily adapt Proposition 4, Lemma 7,
Lemma 8, Propositions 9 and 10 from \cite{MR}:
\begin{prop}\label{gonabla}
Let $(V,\rho)$ be a representation of $\op{GL}(n,q,\mathbb{R})$
and $\rho'$ the induced action on $\mathbb{R}^{n*}\otimes V$. If
$f$ belongs to
$C^{\infty}(P_{\mathcal{F}},V)_{\op{GL}(n,q,\mathbb{R})}$, then
  $\nabla^{\omega_{\mathcal{F}}}f\in C^{\infty}(P_{\mathcal{F}},\R^{n*}\otimes V)_{\op{GL}(n,q,\mathbb{R})}$.
\end{prop}
\begin{proof}
The result is a consequence of the Ad-invariance of the
  Cartan connection $\omega_{\mathcal{F}}$. Indeed :
$$(\nabla^{\omega_{\mathcal{F}}}f)(ug)=\rho'(g)^{-1}(\nabla^{\omega_{\mathcal{F}}}f)(u)\;\forall u\in P_{\mathcal{F}},
  \forall g\in\op{GL}(n,q,\mathbb{R})$$
$$\Longleftrightarrow$$
$$(\nabla^{\omega_{\mathcal{F}}}f)(ug)(X)=[\rho'(g)^{-1}(\nabla^{\omega_{\mathcal{F}}}f)(u)](X)\;\forall
  u\in P_{\mathcal{F}},\forall g\in\op{GL}(n,q,\mathbb{R}),\forall X\in\mathbb{R}^{n}$$
$$\Longleftrightarrow$$
$$(L_{\omega_{\mathcal{F}}^{-1}(X)}f)(ug)=\rho(g^{-1})(L_{\omega_{\mathcal{F}}^{-1}(gX)}f)(u)\;\forall
  u\in P_{\mathcal{F}},\forall g\in\op{GL}(n,q,\mathbb{R}),\forall X\in\mathbb{R}^{n}.$$
If one denotes by $\varphi_{t}$ the flow of
$\omega_{\mathcal{F}}^{-1}(X)$ and by $\varphi_{t}'$ the flow of
$\omega_{\mathcal{F}}^{-1}(gX)$, it suffices then to verify that
$$\frac {d}{dt} f(\varphi_{t}(ug))|_{t=0}=\rho(g^{-1})\frac {d}{dt}
  f(\varphi_{t}'(u))|_{t=0}\;\forall u\in P_{\mathcal{F}},\forall g\in\op{GL}(n,q,\mathbb{R}),$$
or that
$$\varphi_{t}(ug)=\varphi_{t}'(u)g\;\forall u\in P_{\mathcal{F}},\forall g\in\op{GL}(n,q,\mathbb{R}).$$
This property is satisfied : indeed, the fields
$\omega_{\mathcal{F}}^{-1}(gX)$ and $\omega_{\mathcal{F}}^{-1}(X)$
are $R_{g}$-linked because of the $\op{Ad}$-invariance of
$\omega_{\mathcal{F}}$.
\end{proof}
In the same way, we have the following result :
\begin{prop}\label{gonabladiv}
Let $\rho$ be the action of $\op{GL}(q,\mathbb{R})$ on
$S^k(\mathbb{R}^q)$ and $\rho'$ the induced action on
$\mathbb{R}^{q*}\otimes S^k(\mathbb{R}^q)$. If $S\in
C^{\infty}(P_{\mathcal{F}},S^k(\mathbb{R}^n))$ is such that
$(p_{n,q}S)(ug)=\rho(p_{n,q}g^{-1})(p_{n,q}S(u))\;\forall
g\in\op{GL}(n,q,\mathbb{R})$, then
$$(p_{n,q}\nabla^{\omega_{\mathcal{F}}}S)(ug)=\rho'(p_{n,q}g^{-1})(p_{n,q}\nabla^{\omega_{\mathcal{F}}}S(u)).$$
\end{prop}
\begin{proof}
The proof is analogous to the proof of the previous result.
$$(p_{n,q}\nabla^{\omega_{\mathcal{F}}}S)(ug)=\rho'(p_{n,q}g)^{-1}(p_{n,q}\nabla^{\omega_{\mathcal{F}}}S)(u)\;\forall u\in P_{\mathcal{F}},
  \forall g\in\op{GL}(n,q,\mathbb{R})$$
$$\Longleftrightarrow$$
$$(p_{n,q}\nabla^{\omega_{\mathcal{F}}}S)(ug)(X)=[\rho'(p_{n,q}g)^{-1}(p_{n,q}\nabla^{\omega_{\mathcal{F}}}S)(u)](X)\;\forall
  u\in P_{\mathcal{F}},\forall g\in\op{GL}(n,q,\mathbb{R}),\forall X\in\mathbb{R}^{q}$$
$$\Longleftrightarrow$$
$$(L_{\omega_{\mathcal{F}}^{-1}(0,X)}p_{n,q}S)(ug)=\rho(p_{n,q}g^{-1})(L_{\omega_{\mathcal{F}}^{-1}(0,(p_{n,q}g)X)}p_{n,q}S)(u)\;\forall
  u\in P_{\mathcal{F}},\forall g\in\op{GL}(n,q,\mathbb{R}),\forall X\in\mathbb{R}^{q}.$$
If one denotes by $\varphi_{t}$ the flow of
$\omega_{\mathcal{F}}^{-1}(0,X)$ and by $\varphi_{t}'$ the flow of
$\omega_{\mathcal{F}}^{-1}(0,(p_{n,q}g)X)$, it suffices then to
verify that
$$\frac {d}{dt} p_{n,q}S(\varphi_{t}(ug))|_{t=0}=\rho(p_{n,q}g^{-1})\frac {d}{dt}
  p_{n,q}S(\varphi_{t}'(u))|_{t=0}\;\forall u\in P_{\mathcal{F}},\forall g\in\op{GL}(n,q,\mathbb{R}).$$
One concludes using the fact that
$$\varphi_{t}(ug)=\varphi_{t}'(u)g'$$
with $p_{n,q}g'=p_{n,q}g$ because the fields
$\omega_{\mathcal{F}}^{-1}(0,(p_{n,q}g)X)$ and
$\omega_{\mathcal{F}}^{-1}(0,X)$ are $R_{g'}$-linked by $g'$ such
that $p_{n,q}g'=p_{n,q}g$.
\end{proof}
\begin{prop}\label{goinv}
Let $\rho$ be the action of $\op{GL}(q,\mathbb{R})$ on
$S^k(\mathbb{R}^q)$ and $\rho'$ the action on
$S^{k-1}(\mathbb{R}^q)$. If $S\in
C^{\infty}(P_{\mathcal{F}},S^k(\mathbb{R}^n))$ is such that
$(p_{n,q}S)(ug)=\rho(p_{n,q}g^{-1})(p_{n,q}S(u))\;\forall
g\in\op{GL}(n,q,\mathbb{R})$, then
$$(p_{n,q}\op{Div}^{\omega_{\mathcal{F}}}S)(ug)=\rho'(p_{n,q}g^{-1})(p_{n,q}\op{Div}^{\omega_{\mathcal{F}}}S(u)).$$
\end{prop}
\begin{proof}
This can be checked directly from the definition of the divergence
and from the proposition \ref{gonabladiv}. We have successively :
\begin{eqnarray*}
(p_{n,q}\op{Div}^{\omega_{\mathcal{F}}}S)(ug) & = &
\sum_{j=p+1}^{n}(p_{n,q}\nabla^{\omega_{\mathcal{F}}}S)(ug)(p_{n,q}e_{j})(p_{n,q}\epsilon^{j})
\\
& = &
\rho'(p_{n,q}g^{-1})\sum_{j=p+1}^{n}(p_{n,q}\nabla^{\omega_{\mathcal{F}}}S)(u)((p_{n,q}g)p_{n,q}e_{j})(p_{n,q}\epsilon^{j}(p_{n,q}g^{-1}))
\\
                                         & = &
                                         \rho'(p_{n,q}g^{-1})\sum_{j=p+1}^{n}\sum_{i=p+1}^{n}\sum_{r=p+1}^{n}(p_{n,q}\nabla^{\omega_{\mathcal{F}}}S)(u)(p_{n,q}e_{i})(p_{n,q}\epsilon^{r})g_{j}^{i}g_{r}^{-1j}
                                         \\
                                         & = & \rho'(p_{n,q}g^{-1})p_{n,q}\sum_{i=p+1}^{n}(\nabla^{\omega_{\mathcal{F}}}S)(u)(e_{i})(\epsilon^{i}).\\
\end{eqnarray*}
\end{proof}
In our computations, we will make use of the infinitesimal version
of the equivariance relation : if $(V,\tilde{\rho})$ is a
representation of $H(n+1,q+1,\mathbb{R})$ and if $f\in
C^{\infty}(P_{\mathcal{F}},V)_{H(n+1,q+1,\mathbb{R})}$ then one
has
\begin{equation}\label{Invalg}
L_{h^*}f(u) + \tilde{\rho}_*(h)f(u)=0,\quad\forall
h\in\op{gl}(n,q,\R)\oplus\R^{q*}\subset\op{sl}(n+1,\R), \forall
u\in P_{\mathcal{F}}.
\end{equation}
\begin{prop}\label{div1}
For every $S\in C^{\infty}(P_{\mathcal{F}},S^k(\mathbb{R}^n))$
such that
$(p_{n,q}S)(ug)=\rho(p_{n,q}g^{-1})(p_{n,q}S(u))\;\forall
g\in\op{GL}(n,q,\mathbb{R})$, we have
\[p_{n,q}[L_{h^*}\op{Div}^{\omega_{\mathcal{F}}}S - \op{Div}^{\omega_{\mathcal{F}}}L_{h^*}S]=(q+2k-1)p_{n,q}i(0,h)S,\]
for every $h\in\mathbb{R}^{q*}$.
\end{prop}
\begin{proof}
First we remark that the Lie derivative with respect to a vector
field commutes with the evaluation : if
$\eta^1,\ldots,\eta^{k-1}\in \R^{q*}$, we have
\[\begin{array}{lll}(L_{h^*}p_{n,q}\op{Div}^{\omega_{\mathcal{F}}}S)(\eta^1,\ldots,\eta^{k-1})
& =& L_{h^*}(p_{n,q}\op{Div}^{\omega_{\mathcal{F}}}S(\eta^1,\ldots,\eta^{k-1}))\\
 & = &
\sum_{j=p+1}^n(L_{h^*}L_{\omega_{\mathcal{F}}^{-1}(e_j)}p_{n,q}S(p_{n,q}\epsilon^j,\eta^1,\ldots,\eta^{k-1})).\end{array}\]
Now, the definition of a Cartan connection implies the relation
 \[[h^*,\omega_{\mathcal{F}}^{-1}(X)] = \omega_{\mathcal{F}}^{-1}([h,X]),\quad\forall h\in\op{gl}(n,q,\R)\oplus \R^{q*},X\in\R^n,\] where the bracket on the
right is the one of $\op{sl}(n+1,\R)$. It follows that the
expression we have to compute is equal to
\[\sum_{j=p+1}^n(L_{\omega_{\mathcal{F}}^{-1}(e_j)}L_{h^*}p_{n,q}S(p_{n,q}\epsilon^j,\eta^1,\ldots,\eta^{k-1})
+
(L_{[h,e_j]^*}p_{n,q}S)(p_{n,q}\epsilon^j,\eta^1,\ldots,\eta^{k-1})).\]
Finally, we obtain
\[\begin{array}{lll}&& p_{n,q}\op{Div}^{\omega_{\mathcal{F}}}(L_{h^*}S)(\eta^1,\ldots,\eta^{k-1}) -(\rho_{*}(p_{n,q}[h,e_j])p_{n,q}S)(p_{n,q}\epsilon^j,\eta^1,
\ldots,\eta^{k-1})\\
&=&p_{n,q}\op{Div}^{\omega_{\mathcal{F}}}(L_{h^*}S)(\eta^1,\ldots,\eta^{k-1})
+(\rho_{*}(p_{n,q}(h\otimes e_j + \langle h,e_j\rangle
Id))p_{n,q}S)(p_{n,q}\epsilon^j,\eta^1,
\ldots,\eta^{k-1}).\end{array}\] The result then easily follows
from the definition of $\rho$ on $S^k(\R^q)$.
\end{proof}
\begin{prop}\label{div2}
If $S$ is an equivariant function on $P_{\mathcal{F}}$
representing an adapted symbol, we have
\[p_{n,q}[L_{h^*} (\op{Div}^{\omega_{\mathcal{F}}})^l S - (\op{Div}^{\omega_{\mathcal{F}}})^lL_{h^*}S] =
l(q+2k-l)p_{n,q}[i(h) (\op{Div}^{\omega_{\mathcal{F}}})^{l-1}S],\]
for every $h\in \mathbb{R}^{q*}$.
\end{prop}
\begin{proof}
For $l=1$, this is simply the proposition \ref{div1}. Then the
result follows by induction,
 using propositions \ref{goinv} and \ref{div1}. One has indeed, if one supposes the result true to $l-1$,
 that
$$L_{h^*}p_{n,q}\op{Div}^{\omega_{\mathcal{F}}}(\op{Div}^{\omega_{\mathcal{F}}l-1})S-p_{n,q}\op{Div}^{\omega_{\mathcal{F}}l}L_{h^{*}}S$$
is equal to
$$(q+2(k-l+1)-1)p_{n,q}i(h)\op{Div}^{\omega_{\mathcal{F}}l-1}S+p_{n,q}\op{Div}^{\omega_{\mathcal{F}}} L_{h^*}\op{Div}^{\omega_{\mathcal{F}}l-1}S-p_{n,q}\op{Div}^{\omega_{\mathcal{F}}l}L_{h^{*}}S,$$
i.e. to
$$l(q+2k-l)p_{n,q}i(h)\op{Div}^{\omega_{\mathcal{F}}l-1}S.$$
\end{proof}
\begin{prop}\label{nablag1}
If
$f\in\mathcal{C}^{\infty}(P_{\mathcal{F}},\mathbb{R})_{\op{GL}(n,q,\mathbb{R})}$,
then
\[L_{h^{*}}(\nabla^{\omega_{\mathcal{F}}})^{k}f -
(\nabla^{\omega_{\mathcal{F}}})^{k}L_{h^{*}}f
 =-k(k-1)(\nabla^{\omega_{\mathcal{F}}})^{k-1}f\vee
h,\] for every $h\in \mathbb{R}^{q*}$.
\end{prop}
\begin{proof}
If $k=0$, then the formula is obviously true. Then we proceed
 by induction. In view of the symmetry of the expressions that we have to
 compare, it is sufficient to check that they coincide when evaluated on the
$k$-tuple $(X,\ldots,X)$ for every $X\in\R^n$. The proof is
similar to the one
 of proposition \ref{div1} : first the evaluation and the Lie derivative commute :
\[(L_{h^{*}}(\nabla^{\omega_{\mathcal{F}}})^{k}f)(X,\ldots,X)
= L_{h^*}((\nabla^{\omega_{\mathcal{F}}})^{k}f(X,\ldots,X)).\]
Next, we use the definition of the iterated invariant differential
and we let the operators $L_{h^*}$ and
$L_{\omega_{\mathcal{F}}^{-1}(X)}$ commute so that the latter
expression becomes
\[L_{\omega_{\mathcal{F}}^{-1}(X)}L_{h^{*}}((\nabla^{\omega_{\mathcal{F}}})^{k-1}f)(X,\ldots,X)
+(L_{[h,X]^{*}}((\nabla^{\omega_{\mathcal{F}}})^{k-1}f))(X,\ldots,X).\]
By the induction, the first term is equal to
\[(\nabla^{\omega_{\mathcal{F}}})^{k}L_{h^{*}}f(X,\ldots,X)  -(k-1)(k-2)((\nabla^{\omega_{\mathcal{F}}})^{k-1}f\vee
h)(X,\ldots,X).\] For the second term, we use proposition
\ref{gonabla} and relation (\ref{Invalg}) and we obtain, if one
denotes by $\rho$ the action on $S^{k-1}\mathbb{R}^{n*}$,
\[(\rho_{*}((h\otimes X) +\langle h, X\rangle Id)((\nabla^{\omega_{\mathcal{F}}})^{k-1}f))(X,\ldots,X).\]
The result follows by the definition of $\rho_*$.
\end{proof}
\begin{theo}\label{princa}
In the adapted situation, the formula giving the quantization
$Q_{\mathcal{F}}$ is then the following :
\begin{equation}\label{formula}Q_{\mathcal{F}}(\nabla_{\mathcal{F}}, S)(f) =
p_{\mathcal{F}}^{2*^{-1}}(\sum_{l=0}^k C_{k,l}
\langle\op{Div}^{\omega_{\mathcal{F}}^l}
p_{\mathcal{F}}^{2*}S,\nabla_s^{\omega_{\mathcal{F}}^{k-l}}p_{\mathcal{F}}^{2*}f\rangle)\end{equation}
if
\[C_{k,l} =\frac{(k-1)\cdots(k-l)}{(q+2k-1)\cdots(q+2k-l)}\left(\begin{array}{c}k\\l\end{array}\right),\forall l\geq
1,\quad C_{k,0}=1.\]
\end{theo}
\begin{proof}
The proof goes as in \cite{MR}. First, we have to check that the
formula makes sense : the function
\begin{equation}\sum_{l=0}^k C_{k,l} \langle\op{Div}^{\omega_{\mathcal{F}}^l}
p_{\mathcal{F}}^{2*}S,\nabla_s^{\omega_{\mathcal{F}}^{k-l}}p_{\mathcal{F}}^{2*}f\rangle\end{equation}
has to be $H(n+1,q+1,\mathbb{R})$-equivariant. It is obviously
$\op{GL}(n,q,\mathbb{R})$-equivariant by propositions
\ref{gonabla} and \ref{goinv}. It is then sufficient to check that
it is $\mathbb{R}^{q*}$-equivariant. This follows directly from
propositions \ref{div2} and \ref{nablag1} and from the relation
\begin{equation}C_{k,l}l(q+2k-l)=C_{k,l-1}(k-l+1)(k-l).\end{equation}
Next we see, using the results of \cite[p.47]{Capinv}
 that the principal symbol
of $Q_{\mathcal{F}}(\nabla_{\mathcal{F}}, S)$ is exactly $S$, and
formula (\ref{formula}) defines a quantization, that is
projectively invariant,
 by the definition of $\omega_{\mathcal{F}}$.
Next, the naturality of the quantization defined in this way is
easy to understand : it follows from the naturality of the
association of an adapted projective structure $P_{\mathcal{F}}\to
M$ endowed with an adapted normal Cartan connection
 $\omega_{\mathcal{F}}$ to a class of projectively equivalent torsion-free adapted connections on $M$
  and from the naturality of the lift of the equivariant
 functions on $P_{\mathcal{F}}^1M$ to equivariant functions on $P_{\mathcal{F}}$.
\end{proof}
\subsection{Construction in the foliated situation}
In the foliated situation, one can define the invariant
differentiation in this way :
\begin{prop}
The following definition makes sense : if $f$ is a foliated
function on $P(\mathcal{F})$, then
$$(\nabla^{\omega(\mathcal{F})^{k}}f)(v_{1},\ldots,v_{k})=\frac{1}{k!}\sum_{\nu}
L_{\omega(\mathcal{F})^{-1}(v_{\nu_1})}\circ\ldots\circ
L_{\omega(\mathcal{F})^{-1}(v_{\nu_k})}f(u),$$ where
$\omega(\mathcal{F})^{-1}(v)$ is a vector field such that its
image by $\omega(\mathcal{F})$ is equal to $v$.
\end{prop}
\begin{proof}
One has to show that the definition is independent of the choice
of the vector field. Indeed, two such vector fields differ by a
vector field tangent to $\mathcal{F}_{P^{2}N}$ and one can show
that if $f$ is constant along the leaves of
$\mathcal{F}_{P^{2}N}$, then $L_{\omega(\mathcal{F})^{-1}(v)}f$ is
a foliated function too if $v\in\mathbb{R}^{q}$. Indeed, if $X$ is
tangent to $\mathcal{F}_{P^{2}N}$, then
$L_{X}L_{\omega(\mathcal{F})^{-1}(v)}f=0$. To show that, it
suffices to prove that $L_{[X,\omega(\mathcal{F})^{-1}(v)]}f=0$.

\vspace{0.2cm}One has
$0=d\omega(\mathcal{F})(X,\omega(\mathcal{F})^{-1}(v))=X.v-\omega(\mathcal{F})^{-1}(v).\omega(\mathcal{F})(X)-\omega(\mathcal{F})([X,\omega(\mathcal{F})^{-1}(v)])$.
As the first two terms are equal to 0, the third term vanishes
too. One has then that $[X,\omega(\mathcal{F})^{-1}(v)]$ is
tangent to $\mathcal{F}_{P^{2}N}$ and then
$L_{[X,\omega(\mathcal{F})^{-1}(v)]}f=0$.

\end{proof}

\vspace{0.2cm}In the foliated situation, we define the divergence
operator in this way :
\begin{defi} The \emph{Divergence operator}
with respect to the Cartan connection $\omega(\mathcal{F})$ is
defined by
\[\op{Div}^{\omega(\mathcal{F})} : C^{\infty}(P(\mathcal{F}),S^k(\mathbb{R}^q);\mathcal{F}_{P^{2}N})
\to
C^{\infty}(P(\mathcal{F}),S^{k-1}(\mathbb{R}^q);\mathcal{F}_{P^{2}N})
: S\mapsto \sum_{j=1}^q
i(\epsilon^j)\nabla^{\omega(\mathcal{F})}_{e_j}S.\]
\end{defi}
\vspace{0.2cm}One can then easily adapt the propositions
\ref{gonabla}, \ref{goinv}, \ref{div1}, \ref{div2}, \ref{nablag1}.
The proofs of these propositions are completely similar to the
proofs of the corresponding results in \cite{MR}.
\begin{prop}\label{gonabla'}
If $f$ is a $\op{GL}(q,\mathbb{R})$-equivariant foliated function
on $P(\mathcal{F})$ then
  $\nabla^{\omega(\mathcal{F)}}f$ is $\op{GL}(q,\mathbb{R})$-equivariant too.
\end{prop}
\begin{prop}\label{goinv'}
If $$S\in
C^{\infty}(P(\mathcal{F}),S^k(\mathbb{R}^q);\mathcal{F}_{P^{2}N})_{\op{GL}(q,\mathbb{R})},$$
then $$\op{Div}^{\omega(\mathcal{F})}S\in
C^{\infty}(P(\mathcal{F}),S^{k-1}(\mathbb{R}^q);\mathcal{F}_{P^{2}N})_{\op{GL}(q,\mathbb{R})}.$$
\end{prop}
In our computations, we will make use of the infinitesimal version
of the equivariance relation : if $(V,\tilde{\rho})$ is a
representation of $H(q+1,\mathbb{R})$, if $f\in
C^{\infty}(P(\mathcal{F}),V)_{H(q+1,\mathbb{R})}$ then one has
\begin{equation}\label{Invalg'}
L_{h^*}f(u) + \tilde{\rho}_*(h)f(u)=0,\quad\forall
h\in\op{gl}(q,\R)\oplus\R^{q*}\subset\op{sl}(q+1,\R), \forall u\in
P(\mathcal{F}).
\end{equation}
\begin{prop}\label{div1'}
For every $S\in
C^{\infty}(P(\mathcal{F}),S^k(\mathbb{R}^q);\mathcal{F}_{P^{2}N})_{\op{GL}(q,\mathbb{R})}$
we have
\[L_{h^*}\op{Div}^{\omega(\mathcal{F})}S - \op{Div}^{\omega(\mathcal{F})}L_{h^*}S=(q+2k-1)i(h)S,\]
for every $h\in\mathbb{R}^{q*}$.
\end{prop}

\begin{theo}\label{div2'}
For every $S\in
C^{\infty}(P(\mathcal{F}),S^k(\mathbb{R}^q);\mathcal{F}_{P^{2}N})_{\op{GL}(q,\mathbb{R})}$,
we have
\[L_{h^*} (\op{Div}^{\omega(\mathcal{F})})^lS - (\op{Div}^{\omega(\mathcal{F})})^lL_{h^*}S =
l(q+2k-l)i(h) (\op{Div}^{\omega(\mathcal{F})})^{l-1}S,\] for every
$h\in \mathbb{R}^{q*}$.
\end{theo}
\begin{theo}\label{nablag1'}
If
$f\in\mathcal{C}^{\infty}(P(\mathcal{F}),\mathbb{R};\mathcal{F}_{P^{2}N})_{\op{GL}(q,\mathbb{R})}$,
then
\[L_{h^{*}}(\nabla^{\omega(\mathcal{F})})^{k}f -
(\nabla^{\omega(\mathcal{F})})^{k}L_{h^{*}}f
 =-k(k-1)(\nabla^{\omega(\mathcal{F})})^{k-1}f\vee
h,\] for every $h\in \mathbb{R}^{q*}$.
\end{theo}

\begin{theo}\label{princf}
In the foliated situation, the formula giving the quantization
$Q(\mathcal{F})$ is the following :
$$Q(\mathcal{F})(\nabla(\mathcal{F}), S)(f) =
(p^{2}(\mathcal{F}))^{*^{-1}}(\sum_{l=0}^k C_{k,l}
\langle\op{Div}^{\omega(\mathcal{F})^l}
(p^{2}(\mathcal{F}))^*S,\nabla_s^{\omega(\mathcal{F})^{k-l}}(p^{2}(\mathcal{F}))^*f\rangle)$$
 if
\[C_{k,l} =\frac{(k-1)\cdots(k-l)}{(q+2k-1)\cdots(q+2k-l)}\left(\begin{array}{c}k\\l\end{array}\right),\forall l\geq
1,\quad C_{k,0}=1.\]
\end{theo}

\subsection{Quantization commutes with reduction}

\begin{prop}
If $f$ is an equivariant function on $P_{\mathcal{F}}$
representing a basic function, then
$$(\nabla^{\omega_{\mathcal{F}}^{k}}f)(v_{1},\ldots,v_{k})={^{\cal
F}\!\!{\frak p}^2}^{*}(\nabla^{\omega(\mathcal{F})^{k}}(^{\cal
F}\!\hat{\zp}f))(p_{n,q}v_{1},\ldots,p_{n,q}v_{k}).$$
\end{prop}
\begin{proof}
Indeed, one has first that ${^{\cal F}\!\!{\frak
p}^2}_{*}\omega_{\mathcal{F}}^{-1}(v)$ is equal to
$\omega(\mathcal{F})^{-1}(p_{n,q}v)$ modulo a vector field tangent
to $\mathcal{F}_{P^{2}N}$. By induction, if the proposition is
true to $k-1$, it is true for $k$ :
\begin{eqnarray*}
    (\nabla^{\omega_{\mathcal{F}}^{k}}f)(v_{1},\ldots,v_{k}) & = &
    L_{\omega_{\mathcal{F}}^{-1}(v_{k})}(\nabla^{\omega_{\mathcal{F}}^{k-1}}f)(v_{1},\ldots,v_{k-1})
    \\
                                                         & = &
                                                         L_{\omega_{\mathcal{F}}^{-1}(v_{k})}{^{\cal
F}\!\!{\frak p}^2}^{*}(\nabla^{\omega(\mathcal{F})^{k-1}}(^{\cal F}\!\hat{\zp}f))(p_{n,q}v_{1},\ldots,p_{n,q}v_{k-1})\\
                                                         & = & {^{\cal
F}\!\!{\frak p}^2}^{*}(L_{\omega(\mathcal{F})^{-1}(p_{n,q}v_{k})}\nabla^{\omega(\mathcal{F})^{k-1}}(^{\cal F}\!\hat{\zp}f))(p_{n,q}v_{1},\ldots,p_{n,q}v_{k-1})\\
                                                         & = & {^{\cal
F}\!\!{\frak p}^2}^{*}(\nabla^{\omega(\mathcal{F})^{k}}(^{\cal F}\!\hat{\zp}f))(p_{n,q}v_{1},\ldots,p_{n,q}v_{k}).\\
\end{eqnarray*}
\end{proof}
In an other part,
\begin{prop}
If $S$ is an equivariant function on $P_{\mathcal{F}}$
representing an adapted symbol, then
$$p_{n,q}(\op{Div}^{\omega_{\mathcal{F}}^{l}}S)={^{\cal
F}\!\!{\frak p}^2}^{*}(\op{Div}^{\omega(\mathcal{F})^{l}}(^{\cal
F}\!\hat{\zp}S)).$$
\end{prop}
\begin{proof}
Indeed, by induction, if it is true to $l-1$, it is true for $l$ :
\begin{eqnarray*}
   p_{n,q}(\op{Div}^{\omega_{\mathcal{F}}^{l}}S) & = & p_{n,q}\sum_{j=p+1}^{n}i(\epsilon^{j})L_{\omega_{\mathcal{F}}^{-1}(e_{j})}(\op{Div}^{\omega_{\mathcal{F}}^{l-1}}S)
    \\
                                                         & = & \sum_{j=p+1}^{n}L_{\omega_{\mathcal{F}}^{-1}(e_{j})}p_{n,q}i(\epsilon^{j})(\op{Div}^{\omega_{\mathcal{F}}^{l-1}}S)\\
                                                         & = & \sum_{j=p+1}^{n}L_{\omega_{\mathcal{F}}^{-1}(e_{j})}i(p_{n,q}\epsilon^{j}){^{\cal
F}\!\!{\frak p}^2}^{*}(\op{Div}^{\omega(\mathcal{F})^{l-1}}(^{\cal F}\!\hat{\zp}S))\\
                                                         & = & \sum_{j=p+1}^{n}{^{\cal
F}\!\!{\frak p}^2}^{*}L_{\omega(\mathcal{F})^{-1}(p_{n,q}e_{j})}i(p_{n,q}\epsilon^{j})(\op{Div}^{\omega(\mathcal{F})^{l-1}}(^{\cal F}\!\hat{\zp}S)).\\
\end{eqnarray*}
\end{proof}

\begin{theo}
If $S$ is an equivariant function on $P_{\mathcal{F}}$
representing an adapted symbol and if $f$ is an equivariant
function on $P_{\mathcal{F}}$ representing a basic function, then
$$\langle \op{Div}^{\omega_{\mathcal{F}}^{l}}S,
\nabla^{\omega_{\mathcal{F}}^{k-l}}f\rangle={^{\cal F}\!\!{\frak
p}^2}^{*}\langle \op{Div}^{\omega(\mathcal{F})^{l}}(^{\cal
F}\!\hat{\zp}S), \nabla^{\omega(\mathcal{F})^{k-l}}(^{\cal
F}\!\hat{\zp}f)\rangle$$ if $S$ is of degree $k$. The quantization
commutes then with the reduction :
$$Q_{\mathcal{F}}(\nabla_{\mathcal{F}})(S)(f)=Q(\mathcal{F})(\nabla(\mathcal{F}))(^{\cal
F}\!\zp S)(^{\cal F}\!\zp f).$$
\end{theo}

\begin{proof}
Indeed, if
$\op{Div}^{\omega_{\mathcal{F}}^{l}}S=v_{1}\vee\ldots\vee
v_{k-l}$,
\begin{eqnarray*}
  \langle\op{Div}^{\omega_{\mathcal{F}}^{l}}S,
\nabla^{\omega_{\mathcal{F}}^{k-l}}f\rangle  & = &
(\nabla^{\omega_{\mathcal{F}}^{k-l}}f)(v_1,\ldots,v_{k-l})
    \\
                                                         & = & {^{\cal
F}\!\!{\frak p}^2}^{*}(\nabla^{\omega(\mathcal{F})^{k-l}}(^{\cal F}\!\hat{\zp}f))(p_{n,q}v_1,\ldots,p_{n,q}v_{k-l})\\
                                                         & = & \langle p_{n,q}\op{Div}^{\omega_{\mathcal{F}}^{l}}S,
{^{\cal
F}\!\!{\frak p}^2}^{*}\nabla^{\omega(\mathcal{F})^{k-l}}(^{\cal F}\!\hat{\zp}f)\rangle.\\
\end{eqnarray*}
The conclusion follows then from Theorems \ref{princa} and
\ref{princf}.
\end{proof}


\begin{thebibliography}{Dillo 83}


\bibitem[BHMP02]{BHMP} Boniver F, Hansoul S, Mathonet P, Poncin N,
{\it Equivariant symbol calculus for differential operators acting
on forms}, Lett. Math. Phys., {\bf 62}(3) (2002), pp 219-232

\bibitem[BM01]{BM} Boniver F, Mathonet P, {\it Maximal subalgebras of vector fields for equivariant
quantizations}, J. Math. Phys., {\bf 42}(2) (2001), pp 582-589

\bibitem[BM06]{BM2} Boniver F, Mathonet P, {\it IFFT-equivariant
quantizations}, J. Geom. Phys., {\bf 56} (2006), pp 712-730

\bibitem[Bor02]{MB} Bordemann M, {\it Sur l'existence d'une prescription d'ordre naturelle projectivement
invariante}, arXiv:math.DG/0208171

\bibitem[BHP06]{BHPSingReduDefoQuan} Bordemann M, Herbig H-C,
Pflaum M, {\it A homological approach to singular reduction in
deformation quantization}, Singularity theory World Sci. Publ.,
Hackensack, NJ (2007), pp 443--461

\bibitem[Blum84]{Blum} Blumenthal R.A., {\it Cartan connections in
foliated bundles}, Michigan Math. J., {\bf 31} (1984), pp 55-63

\bibitem[CSS97]{Capinv} $\mathrm{\check{C}}$ap A, Slov\'{a}k J, Sou$\mathrm{\check{c}}$ek V, {\it Invariant operators on manifolds
  with almost Hermitian symmetric structures. I. Invariant differentiation.}, Acta
  Math. Univ. Comenian. (N.S.), {\bf 66}(1) (1997), pp 33-69

\bibitem[DLO99]{DLO} Duval C, Lecomte P, Ovsienko V, {\it Conformally equivariant quantization: existence and
uniqueness}, Ann. Inst. Fourier, {\bf 49}(6) (1999), pp 1999-2029

\bibitem[DO01]{DO} Duval C, Ovsienko V, {\it Projectively equivariant quantization and symbol calculus: noncommutative hypergeometric
functions}, Lett. Math. Phys., {\bf 57}(1) (2001), pp 61-67

\bibitem[GP07]{GP} Grabowski J, Poncin N, {\it On quantum and classical Poisson algebras}, Banach Center Publ.,
{\bf76} (2007); activity SMR 1665, ICTP, Trieste
(http://poisson.zetamu.com)

\bibitem[Han06]{SH} Hansoul S, {\it Existence of natural and projectively equivariant
quantizations}, Adv. Math., {\bf 214}(2) (2007), pp 832--864

\bibitem[Hui07]{HLSingQuanCommRedu} Hui L, {\it Singular unitary in
``quantization commutes with reduction''}, ArXiv: 0706.1471

\bibitem[Hue02]{JHQuantReduc} Huebschmann J, {\it Quantization and
Reduction}, ArXiv: math/0207166

\bibitem[Hue06]{JHSingQuant} Huebschmann J, {\it Classical phase space singularities and
quantization}, ArXiv: math-ph/0610047v1

\bibitem[HRS07]{JRSQuantStrat} Huebschmann J, Rudolph G, Schmidt M, {\it A gauge model for quantum mechanics on a stratified
space}, ArXiv: hep-th/0702017

\bibitem[Koba72]{Koba} Kobayashi S, {\it Transformation groups in differential
  geometry}, Springer-Verlag, New York (1972), Ergebnisse der Mathematik und
  ihrer Grenzgebiete, Band 70.

\bibitem[Lec00]{PL} Lecomte P, {\it On the cohomology of $\op{sl}(m+1,\R)$ acting on differential operators and $\op{sl}(m+1,\R)$-equivariant
symbol}, Indag. Math., {\bf 11}(1) (2000), pp 95-114

\bibitem[Lec01]{PL2} Lecomte P, {\it Towards projectively equivariant
quantization}, Progr. Theoret. Phys. Suppl., {\bf 144} (2001), pp
125-132

\bibitem[LMT96]{LMT} Lecomte P, Mathonet P, Tousset E, {\it Comparison of some modules of the Lie algebra of vector fields}, Indag. Math., {\bf 7}(4) (1996), pp 461-471

\bibitem[LO99]{LO} Lecomte P, Ovsienko V, {\it Projectively equivariant symbol
calculus}, Lett. Math. Phys., {\bf 49}(3)(1999), pp 173-196

\bibitem[MR05]{MR} Mathonet P, Radoux F, \textit{Natural and projectively equivariant quantizations by means of Cartan
connections}, Lett. Math. Phys., {\bf 72}(3) (2005), pp 183-196

%\bibitem[OO93]{OO} Ovsienko V Yu, Ovsienko O D, {\it Projective structures and infinite-dimensional Lie algebras associated with a
%contact manifold}, Advances in Soviet Mathematics, {\bf 17} (1993)

%\bibitem[OR92]{OR} Ovsienko V, Roger C, \textit{Deformations of Poisson brackets and extensions of Lie algebras of contact vector
%fields}, Uspekhi Mat. Nauk, {\bf 47}(6) (1992), pp 141-194,
%Russian Math. Surveys {\bf 47}(6) (1992), pp 135-191

%\bibitem[Ovs05]{VO} Ovsienko V, \textit{Vector fields in presence of a contact
%structure}, Enseign. Math. (2){\bf 52}(3-4) (2006), pp 215--229

\bibitem[Pfl02]{PflDefoQuanSympOrbi} Pflaum M J, {\it On the deformation quantization of symplectic
orbispaces}, Differential Geom. Appl. {\bf 19}(3) (2003), pp
343--368

\bibitem[Rich01]{RichTransGeom} Richardson K, {\it The transverse geometry
of G-manifolds and Riemannian foliations}, Illinois J. Math., {\bf
45}(2) (2001), 517-535.

\bibitem[Wol89]{RW} Wolak R, {\it Foliated and associated geometric structures on foliated manifolds}, Ann. Fac. Sci. Toulouse {\bf 10} (3)
(1989), pp 337-360
\end{thebibliography}
\end{document}